\documentclass{article}

\usepackage{color}
\usepackage{graphicx}
\usepackage[normalem]{ulem}
\usepackage{morefloats}
\usepackage{hyperref}
\hypersetup{
    colorlinks=true,
    linkcolor=blue,
    filecolor=magenta,      
    urlcolor=cyan,
}
\newtheorem{conjecture}{Conjecture}
\author{Barry Brent}

\date{15h 29 March 2017, preliminary draft}
\begin{document}

\title{Experiments with the dynamics of the Riemann zeta function}

\maketitle
\begin{abstract}
We collect experimental evidence 
for several propositions,
including the following:
(1) For each Riemann zero $\rho$
(trivial or nontrivial)
and each zeta fixed point
$\psi$ 
there is a 
nearly logarithmic
spiral $s_{\rho, \psi}$
with center $\psi$
containing $\rho$. 
(2) $s_{\rho, \psi}$ interpolates a 
 subset  $B_{\rho, \psi}$ of the
backward zeta orbit of $\rho$
comprising a set
of zeros of all iterates of  zeta.
(3) If zeta is viewed as a function on sets, 
$\zeta(B_{\rho, \psi}) = B_{\rho, \psi} \cup \{ 0 \}$.
(4) $B_{\rho, \psi}$ has nearly
uniform angular
distribution around 
the center of $s_{\rho, \psi}$.
We will make these statements precise.
\end{abstract}

\bibliographystyle{plain}

\section{\sc introduction} 
\subsection{Overview.}
 For complex $w$, let $\zeta^{\circ -}(w)$
denote the backward zeta orbit
of $w$  so that
$\zeta^{\circ -}(w) = \{s \in \bf{C}$  s.t. 
some iterate of zeta takes $s$ to $w \}$.
If the sequence
$B =  (a_0,  a_1, a_2,  ... )$
satisfies $a_0 = w$
and $\zeta(a_n) = a_{n-1}$
for all $n  \geq 1$, we
say that  $B$
is a  branch of $\zeta^{\circ -}(w)$.
If $B$ converges with  $\lim B =  \lambda$, say,
we conjecture that $B$ is unique and we write 
$B$ as $B_{w, \lambda}$.
We collect numerical evidence 
supporting the  following claims, which
will be made precise below:
(1) that for each of a countable set
of  non-real
zeta fixed points $\psi$ and each 
nontrivial Riemann zero
$\rho$,
$B_{\rho, \psi}$ exists, is unique, and
is the center of a nearly logarithmic
spiral,
say  $s_{\rho, \psi}$,
interpolating  $B_{\rho, \psi}$;
(2) that  the members of
 $B_{\rho, \psi}$
 are distributed nearly uniformly on 
  $s_{\rho, \psi}$;
  (3) and that there is another set of
  real zeta fixed points $\psi = \psi_{-2n}$
  near the trivial zeros $-2n$, $-2n \leq -20$,
  such that consecutive members of 
  $B_{\rho, \psi_{-2n}}$
  rotate around  $\psi_{-2n}$ 
  through an angle of $\approx \pi$ or $2 \pi$,
  depending upon the parity of $n$,
 so that the members of 
  $B_{\rho, \psi_{-2n}}$ lie on
  a curve that  is very nearly a straight line
  passing through  both $\psi_{-2n}$
  and $\rho$.
  \newline \newline
  We  treat in detail 
 relationships between
  zeta basins of attraction  and 
  the branches of $\zeta^{\circ}(\rho)$
  that we have  observed experimentally.
  The resulting plots are  included here because
  they suggested the  presence of the spirals
  which are the main subject of the article,  and
  so--we speculate--they may eventually also suggest
  the ideas needed to analyze these spirals;
for we have not proved any theorems.
\newline \newline
Our experiments were done with \it Mathematica \rm and 
spot-checked with 
\it Sage\rm.  Data files and \it Mathematica  \rm notebooks are 
posted on the ResearchGate site \cite{Br}.
\subsection{Possible bearing on the Riemann hypothesis.}
Here are several scenarios.
(1)  A clear understanding  of the spirals $s_{\rho, \psi}$
with 
zeta fixed point centers $\psi$ and 
passing through Riemann zeros 
$\rho$ might  lead to a sort of dictionary, so that
the Riemann hypothesis might be 
put in a form that speaks of zeta fixed points instead of 
zeta zeros. (2)   If  (as we conjecture)
the spirals $s_{\rho, \psi}$
are  approximated by logarithmic 
spirals, it might be possible to
confine the $s_{\rho, \psi}$ to
lie within 
spiral-shaped ``error bands'' about 
logarithmic spirals.
Because each zero $\rho$ lies on each 
$s_{\rho, \psi}$,
each zero would lie in each one of 
a countable collection of these
error bands in the complex plane
(one for each zeta fixed point $\psi$);
then the zeros would
be confined to
the intersection of these error
bands,  and this region 
might be small. Our very incomplete
knowledge
of the $s_{\rho, \psi}$  
is founded on prior knowledge of
the  locations of the zeros, so that
this scheme is tainted with circularity;
but perhaps this taint might in
some way
be removed. (3) 
Under the conjectures  stated in this
article, measuring from the fixed points $\psi$,
a zero $\rho$ lies at the ``first'' intersection
(in terms of arc length, say)
of  $s_{\rho, \psi}$ with the critical line;
and so the Riemann hypothesis
might be restated  in terms of
these intersections. 
Each $\rho$ appears to lie on
all  of the $s_{\rho, \psi}$,
and so we might  eventually
obtain a countable family of
conditions on these intersections,
which could, possibly, be 
in some way usefully 
combined. Section 6.1.3 discusses 
some data that support 
Conjecture 4, which  codifies a part of
this scenario.
\subsection{More definitions.}
Let $\zeta$ denote the Riemann zeta function
and let us write the iterates of
a function $f$ as 
$f^{\circ 0}(z) = z$ and 
 $f^{\circ (n+1)}(z)= f(f^{\circ n}(z))$
for $n =0, 1, ....$
An $n$-cycle for $f$ is an
$n$-tuple $(c_0, ... , c_{n-1})$
such that $f(c_{n-1}) = c_0$ and 
$f(c_k) = c_{k + 1}$
when  $k \neq n - 1.$
The forward orbit  of $w$ under $f$ 
is the sequence
$(w, f(w), f^{\circ 2}(w), ... )$.
The backward orbit of $w$ under $f$ 
is the set 
of complex numbers $s$ such that
$f^{\circ n}(s) = w$ for some
integer $n \geq 0$. 
Let the symbol $f^{\circ -}(w)$
 denote this 
backward orbit; if
$w$ does not belong to a cycle, 
$f^{\circ -}(w)$ carries the 
 structure of a rooted tree in 
 which the root is $w$ and the children of 
$s \in  f^{\circ -}(w)$
are the solutions $t$ of  $f(t) = s$.
We will call any path in 
$f^{\circ -}(w)$
that begins at $w$ a branch of
$f^{\circ -}(w)$ 
(also: ``a branch of the inverse of $f$''.)
Such a branch, then, is a sequence
$(a_0, a_1, a_2, ....)$ with
$a_0 = w$ and $a_n = f(a_{n+1})$ for
each non-negative integer $n$.
Since the Riemann hypothesis has been verified 
within the range of our observations,
we write without
ambiguity $\rho_k$ for
the $k^{th}$ nontrivial Riemann zero 
ordered by height
above the real axis and 
$\rho_{-k}$ for its complex conjugate. 
If  $\zeta^{\circ n}(z) = 0$ 
and  $\zeta^{\circ n-1}(z)$ is a nontrivial
Riemann zero, we call $z$ a
nontrivial zero of $\zeta^{\circ n}$.
For typographical reasons,
we will occasionally write 
 $\zeta_n$ for $\zeta^{\circ n}$;
within our article, there should no confusion 
with other common uses
of this symbol.
\newline \newline
For $z \in \bf{C} \rm \cup \{\infty\}$, 
$A_z :=\{u \in \bf{C}$ s.t.
$\lim_{n \to\infty} 
\zeta^{\circ n}(u) = z\}$
(the ``basin of attraction'' of $z$ under 
zeta iteration.)
Let $\phi \approx -.295905$
be the largest negative
zeta fixed point.
Then $A_{\phi}$ is a fractal \cite{W}; 
each nontrivial Riemann zero 
appears to lie in an irregularly
shaped bulb of $A_{\phi}$ 
(Figures 1.1, 3.2, and section 3 more generally.)
\newline \newline
For a spiral
$s$  with center $\gamma$,
let  $\alpha_s$ be  the point on
the  intersection of the critical line
and $s$ closest to $\gamma$. 
We  define a real-valued function
$\theta(z)$
on complex numbers $z \in s, 
|z-\gamma| \leq |\alpha_s - \gamma|$ 
by requiring that 
$\theta(z)  
\equiv \arg(z - \gamma) \pmod{2\pi}$,
and that $\theta(z)$
increases continuously and monotonically as $z$
moves around $s$ from $\alpha$ in the direction 
of decreasing $|z-\gamma|$.
In other words,
$\theta(z)$ behaves up to a multiplicative 
constant like a  winding number.
If  a sequence  $(a_1,  a_2, ...)$
lies on  a spiral 
$s$ with center $\gamma$
and for some pair of real numbers $A > 0, B >0$
and all 
 $k = 1, 2, ..., |\theta(a_k) - 
\theta(a_{k+1})| < A e^{-Bk}$,
we will say that the  $a_k$
are distributed nearly uniformly around
$s$.
\newline \newline
For complex $z$,
let $r(z) := |z-\gamma|$. Let $m, b$ be real numbers,
so that $r(z) = \exp (m\theta(z) + b)$ 
describes a logarithmic spiral 
with center
$\gamma$ and typical element
$\exp (m\theta(z) + b) \exp (i \theta(z))$. 
Suppressing the dependence on $m$ and $b$, let
$$d_{rel}(\gamma, z) := 
\left | \frac{z-\gamma-\exp (m\theta(z) + b) \exp (i \theta(z))}{z-\gamma}  \right|.
$$
We  say that $s$ is  $c$-nearly logarithmic  for
real positive $c$ if
$$\max_{z \in s, 0 < |z - \gamma| \leq |\alpha_s - \gamma|} d(\gamma, z) < c$$
for some  $m$ and $b$.
\newline \newline
We require a notion of  ``very nearly a straight line.''
Suppose (1) a complex curve $C$ of finite arc length
has an initial point $z_I$  and terminates at  a point   $z_T$,
(2) that there are real numbers $m$ and $b$
such that 
$$
\lim_{z \in C, z \to z_T} 
 \left|\frac{ \Im(z) - (m \Re(z) + b)}{\Im(z)}\right| = 0,
$$
and (3) that the convergence has exponential decay as 
$|z - z_T|$ decreases from $|z_I - z_T|$ to zero.  
Then we say that  $C$ is very nearly a straight line.
\newline \newline
We need measures of the absolute
and relative deviations
of points $a_k$ in a branch $B_{\rho,\psi}$ 
of the backward orbit of a
nontrivial Riemann zero  $\rho$ from
a logarithmic spiral  fitted to that branch
using \it Mathematica\rm's FindFit command.
Suppose the $a_k$ are interpolated by a spiral
$s_{\rho, \psi}$ centered at $\psi$
such that for 
$z \in s_{\rho, \psi},
r(z) = |z  - \psi|$
and $\theta(z) = \arg (z - \psi)$.
Further, suppose that $s_{\rho, \psi}$
has a log-linear model 
$\tilde{r}(z) = \exp (m\theta(z) + b)$ 
for real numbers $m$ and $b$,
in which we have
fitted the points $(\theta(a_k), \log r(a_k))$
to a straight line.
Then we will write
$$d_{abs}(\rho,\psi,k):= |a_k  - \psi-\tilde{r}(a_k) e^{i \theta(a_k)}|$$
and
$$d_{rel}(\rho,\psi,k):=
\left |\frac {a_k - \psi - 
\tilde{r}(a_k) e^{i \theta(a_k)}}{a_k - \psi}\right |.$$
(This is an  abuse of our earlier notation for
$d_{rel}(\gamma, z)$ which should not be confusing.)
\begin{figure}[!htbp]
\centering
\includegraphics[scale=1]{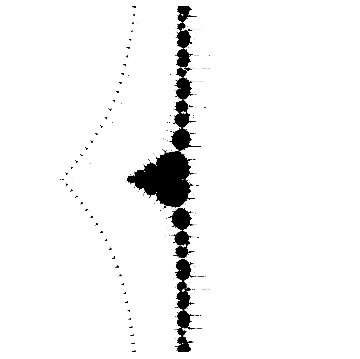}
\vskip .07in
\sc{fig. 1.1: $A_{\phi}$}
\end{figure}
\subsection{Summary of observations.}
We express 
some of our observations as 
explicit conjectures.
\begin{conjecture}
For each pair of positive integers $L$ and $n$, 
there is 
a Riemann zeta $L$-cycle  $\Lambda$
such that  
 the $\lambda \in \Lambda$ with
maximum imaginary part is  close to  the $n^{th}$ nontrivial zero 
$\rho_n$ (but we will not define this use of 
the term ``close''
more precisely in this draft)
and with the following properties.
\newline \newline
(1)  Each $\lambda \in \Lambda$ is 
a repelling fixed point of 
$\zeta^{\circ L}$ lying
in the intersection of boundaries  
$\partial A_{\phi} \cap \partial A_{\infty}$ 
in the usual topology on $\bf{C}$.
\newline \newline
(2) For  each $\lambda \in \Lambda$
there is a  complex number $z_{\lambda}$
and a natural number $0 \leq j_{\lambda} \leq L-1$
such that 
\newline
(a) $\zeta^{\circ j_{\lambda}}(z_{\lambda}) = \rho_n$.
\newline
(b)  For some small positive $c$,
$\lambda$ is the
 center of a 
 $c$-nearly logarithmic spiral
 $s_{z_{\lambda}, \lambda}$
 interpolating a branch $B_{z_{\lambda},\lambda}$
 of $(\zeta^{\circ L})^{\circ -}(z_{\lambda})$
 (but we will not be more precise about this use of the term ``small'' in this draft.)
 \newline
 (c) $\lim B_{z_{\lambda},\lambda} = \lambda$.
 \newline
 (d) $\bigcup_{\lambda \in \Lambda} 
 B_{z_{\lambda},\lambda}$ is a branch  $B_{\rho_n,\Lambda}$ of 
 $\zeta^{\circ -}(\rho_n)$.
 \newline 
 (e) The members of  $B_{z_{\lambda},\lambda}$
 are distributed nearly uniformly on 
  $s_{z_{\lambda}, \lambda}$.
  \newline \newline
(3) If  $z \in \zeta^{\circ -}(\rho_n)$ 
 then for
 some positive integer $j$ and some positive integer $L$
 and some $L$-cycle $\Lambda, \zeta^{\circ j}(z) \in B_{\rho_n,\Lambda}$.
\end{conjecture}
Conjecture 1 is essentially a
description of patterns we observed reliably
in  numerous experiments; clause 2e, in particular,
is  plausible on its face in  view of the spiral plots
exhibited in Figures  5.2,  5.3,  5.4, 5.7, 5.9, and 6.1.  Figure 5.5 supports, in particular,  our use in this clause of the term ``nearly uniform''  as we have
defined it above.
\begin{conjecture}
 (a) If $L = 1$, so that the unique element of
 $\Lambda$ is close to $\rho_n$,
 then we
 write $\Lambda =
 \{\psi_n \}$, and we have that
 $j_{\psi_n} = 0$ and 
 $z_{\psi_n}= \rho_n$.
 In this situation, in clause (b)
 of Conjecture 1 we can take $c< e^{-4}$;
 furthermore, if $L = 1$, 
 then  the infimum of the  valid values of $c$ 
 goes  to zero as $n \to \infty$.
\end{conjecture}
Some evidence for this conjecture appears in section 6.1.1.
\newline \newline
When $L = 1$ we have restricted our claims in
this conjecture to 
spirals  $s_{\rho_n, \psi_n}$ because that is the case
we have checked most thoroughly, but we have also  checked  less thoroughly
the case in which the zero is fixed (usually, $\rho_1$)
and $\psi_n$ varies over positive $n$. It appears to us
that the conjecture generalizes to all pairs 
$\rho =\rho_{n_1}, \psi = \psi_{n_2}, (n_1, n_2) \in \bf{Z}^{\geq 1} \times \bf{Z}^{\geq 1}$.
\newline \newline
Because $\lambda \in \Lambda$ would be repelling, it  
would also be an attracting 
fixed point of a local
branch of the functional inverse
of $\zeta^{\circ L}$, 
and then the convergence
would follow
from standard results, 
for example, Theorem 2.6 of \cite{HY}.
 \begin{conjecture}
(1) There are
repelling zeta fixed points
$\psi_{-2n}$ near the trivial zeros
$-2n \leq -20$. 
\newline
(2) For any  nontrivial Riemann zero 
$\rho$, members of  
$B_{\rho, \psi_{-2n}}$
 lie on a curve which is 
  very nearly a straight line
segment.
If $n$ is even, 
the endpoints are $\rho$ and
$\psi_{-2n}$; otherwise, 
the curve passes through $\psi_{-2n}$
and terminates at $\rho$.
\newline
(3) If $2n \equiv 0$ (mod $4)$,
then $$\arg \frac {d\zeta(\psi_{2n})}{dz}
\approx 2 \pi;$$ 
If $2n \equiv 2$ (mod $4$),
then $$\arg \frac {d\zeta(\psi_{2n})}{dz}
\approx  \pi.$$
\end{conjecture}
(See Figure 5.7. Some evidence for 
this conjecture appears in section 6.1.2.)
These observations are
consistent, of course,
with the proposition that 
 $B_{\rho, \psi_{-2n}}$
 is interpolated by a spiral.
 Our  computations of $\zeta(x)-x$,
 $x$ real, indicate 
 that the $\psi_{-2n}$ are real
 (the graph crosses the $x$-axis near each
 trivial zero we examined.)
It was conceivable that there might
be a broader relationship of the same kind 
between the derivative of zeta at a given
fixed 
point and the structure of the associated spirals 
centered at those fixed points elsewhere
in the complex plane, but our experiments
have not verified any such relationship.
It is suggestive, of course,
that these relationships are
exact when zeta is considered 
as a function of real numbers
and $\psi_{-2n}$ is replaced 
by $-2n$.
\newline \newline
The following conjecture codifies
in part the scenario of section 1.2
\begin{conjecture}
(1) The  relative deviation $d_{rel}(\rho_n,\psi_n,0)$
of  the nontrivial Riemann zero
$\rho_n$ from  a logarithmic spiral
fitted  to  the elements of
$B_{\rho_n,\psi_n} = (a_0, a_1, ...)$
satisfies
$$\lim_{n \to \infty} d_{rel}(\rho_n,\psi_n,0) = 0.$$
In particular, if we write
$$D_{rel}(N) =  \log \frac 1N \sum_{n=1}^N  d_{rel}(\rho_n,\psi_n,0),$$
then there exist two exponents  
$0 < e_1 < e_2 < 1$ such that
$-(\log N)^{e_1} < D_{rel} (N)< -(\log N)^{e_2}$
for $N = 1, 2, ....$
\newline
(2)  The absolute deviation  $d_{abs}(\rho_n,\psi_n,0)$
satisfies
$$\lim_{n \to \infty} d_{abs}(\rho_n,\psi_n,0) = 0.$$
In particular, if we write
$$D_{abs}(N) =  \log \frac 1N \sum_{n=1}^N  d_{abs}(\rho_n,\psi_n,0),$$
then  $\frac 1N <  D_{abs}(N) <  \sqrt{\frac 1N}$.
\end{conjecture}
We  have provided some support for
clause (1)  in Figure 6.5 of section 6.1.3, using $e_1 = .8$
and $e_2 =  .85$;  clause (2) is supported by  the
plot in the right panel of Figure 6.7 in the same section.
\newline \newline
We examined the possibility that
branches of 
$\zeta^{\circ -}(z), \zeta(z) \neq 0$
for arbitrary $z$ on the critical line converge to
the same fixed points as 
branches of  $\zeta^{\circ -}(\rho), \rho$
a non-trivial zero.
We tested various such $z$ 
and found spiral
branches converging to the fixed points of zeta.
So  it ought be
possible to explain the spirals
with a theory that avoids
any appeal to special properties of the 
Riemann zeros. We are not
going to describe these experiments
in any further detail in this article.
\newline  \newline
An analogy from fluid mechanics led us to
check for invariance of  branches of
$\zeta^{\circ -}(z)$
for these $z$
under rotation about 
the fixed points at their centers.
We found that the deviation from
this sort of invariance is
systematic and can itself
be described by referring
to (other) logarithmic
spirals.
\newline \newline
We made a brief survey of functions
 other than zeta to gauge the 
 extent of the spiral phenomenon,
 which we will not describe in
  any further detail
 than the following.
 Functions as simple as cosine appear
 to exhibit this behavior.
 We also observed it in, for example, the
 Ramanujan $L$-function.
 We  hope to carry out another survey
 with a different software package.
\subsection{Prior work.}
Many authors
have examined
the Riemann zeta
function with computers.
Notable citations from the perspective of
this article are
Arias-de-Reyna \cite{A},
Broughan \cite{B},
Cloitre, \cite{C},
Kawahira \cite{Ka}, King \cite{Ki, Ki2}
and Woon \cite{W}.
\section{\sc methods}
\subsection{Quadrant plots.}
We will be displaying
colored plots (say, ``quadrant plots'')  
depicting, for a
point $w$ of $\bf{C}$ and a 
meromorphic function $f$, 
the quadrant of $f(w)$.  We use 
quadrant plots in three ways:
(1) to determine small squares
containing exactly one solution of
an equation of interest, so that
this information can be used by standard
equation-solving routines to find a solution
to several hundred digits of precision (which we find
is necessary, for example, to locate zeta cycles)
lying in a particular region;
(2) to superimpose quadrant plots  upon plots of
the basin of attraction $A_{\phi}$. These
two kinds of plot typically interlock in a way that
helps us to  understand the meaning of
many small irregular features of $A_{\phi}$;
and (3) to  show how the quadrant plots
spiral as we reduce the size of the plot 
window about a fixed point of zeta or one of
its iterates. Observation of
these spiral motions  was our
first indication that forward orbits near 
fixed points do lie on spirals.
\newline \newline
In 
quadrant plots, the 
boundaries of 
single-colored regions are
$f$ pre-images of the axes--curves 
corresponding to
zero sets of $\Re(f(s))$ and $\Im(f(s))$; 
the apparent intersections signal the
presence of zeros or poles of $f$. 
By adjusting the color scheme to 
distinguish between
regions where $|f|$ 
is large or small,
we can try to distinguish zeros
from poles. Some apparent
intersections
are revealed to be illusions by 
a change of scale.
Similar but colorless methods
for plotting zeros were
put to use in \cite{A}.
\newline \newline
The visualized region is partitioned 
into small squares,
each of which is represented by
a pixel. 
We choose a test point $s$
in each square.
The pixel representing the square
is colored according to the rules in 
Table 1.
In the table,
the region $D$ is a disk with center $s = 0$ and 
large radius $r$ (chosen as may be convenient.) 
We denote the complement of $D$ as $-D$.
\newline \newline 
\hskip 1in
\begin{tabular}{|c|c|} \hline
Location of $f(s)$& Color of pixel depicting region containing 
$s$\\ \hline \hline
real and imaginary axes & black \\ \hline 
$D \hskip .05in\cap $ Quadrant I  &  rich blue\\ \hline 
$ - D \hskip .05in\cap$ Quadrant I  &  pale blue\\ \hline
$D \hskip .05in\cap $ Quadrant II  & rich red \\  \hline
$ - D \hskip .05in\cap$ Quadrant II  &  pale red\\ \hline
$D \hskip .05in\cap $ Quadrant III  & rich yellow \\ \hline 
$ - D \hskip .05in\cap$ Quadrant III  &  pale yellow\\ \hline
$D \hskip .05in\cap $ Quadrant IV  &  rich green\\ \hline
$ - D \hskip .05in\cap$ Quadrant IV  &  pale green\\ \hline 
\end{tabular} 
\vskip .1in
\sc{table 1: coloring scheme for quadrant plots}
\rm
\vskip .1in 
\hskip -.19in
The junction of four rich colors 
represents a zero,
the junction of four pale colors represents a pole,
and
the boundary of two appropriately-colored regions is an 
$f$ pre-image of an
axis. An example is shown in Figure 2.1: 
$s \mapsto (s - 1)^2 (s - i) (s + 1)^5/(s+i)^3$.
(We have superimposed a pair of axes on this 
quadrant plot.)
\begin{figure}[!htbp]
\centering
\includegraphics[scale=.8]{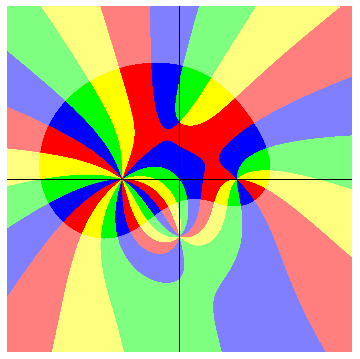}
\vskip .07in
\sc{fig. 2.1: quadrant plot of $s \mapsto 
(s - 1)^2 (s - i) (s + 1)^5/(s+i)^3$}
\end{figure}
\newpage
\subsection{High precision equation solving under geometric constraints.}
\begin{figure}[!htbp]
\centering
\includegraphics[scale=.3]{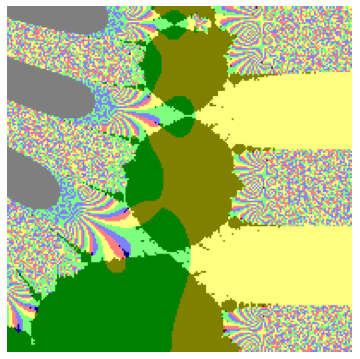}
\hskip .07in
\includegraphics[scale=.3]{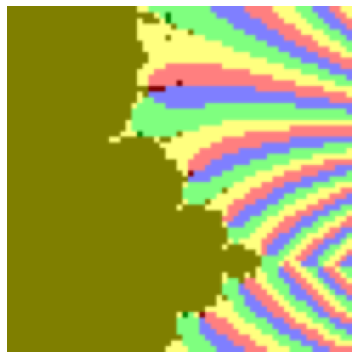}
\hskip .07in
\includegraphics[scale=.3]{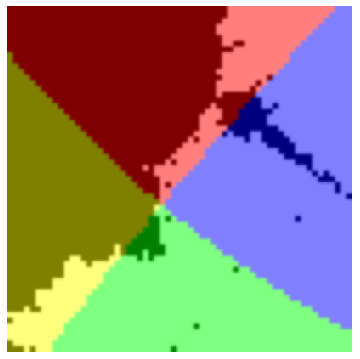}
\hskip .07in
\vskip .1in
\sc{fig. 2.2: zooming in on $\lambda_1$  via quadrant plots of $\zeta^{\circ 3}(s) - s$}
\end{figure}
We illustrate the  application  of quadrant
plots to solve equations under geometric  constraints
by  showing how we found the three-cycle 
$\Lambda$ described in section 5.4. In Figure
2.2, we have superimposed  plots of $A_{\phi}$
and  quadrant plots of 
$s \mapsto \zeta^{\circ 3}(s) - s$
on small squares near the first non-trivial Riemann zero.
The resolution has been kept low to speed up the computations;
high resolution is not particularly helpful in  this
situation. The left  panel is a square with side length $20$
and center $\rho_1$. The central panel depicts a square 
with side length $2$ and center $\rho_1 + 3.4 +.1i$;  we have 
adjusted the center to keep in view a particular four-color
junction visible in the left panel. It represents a solution
of  $\zeta^{\circ 3}(s) - s = 0$,  namely the  three-cycle element 
$\lambda_1$ we are trying to compute. At this stage, if we used,
for example, the \it Mathematica \rm command 
FindRoot constrained to search within this square, 
it might land on any one of the several  four-color 
junctions we see in the central panel.  So we
change the center of the plot again, this time to
$\rho_1 + 3.46 +.103 i$, and  we  make a 
square plot centered there with side length $.02$.
This is shown in the right panel of Figure 2.2.
Next we use a slow, ``handmade'' routine
to approximate $\lambda_1$ by searching within this square.
Then, using this approximation as the beginning value for a 
search with FindRoot, we obtain a solution with $500$ digits
of precision:
\newline
$3.9589623348847434673516458439896123461477039951866801455882506555054$
\newline
$331235719797619129160432526832126428515417856326242408422124490775895$
\newline
$37219976674458409141742662175701089081252727395073714398968532356378$
\newline
$12255138302084634149524670809965144703541657360428502230820135428609$
\newline
$38536894453944241116438492746243199878001238993540770158034816978947$
\newline
$866042863811536518002674033394246742728451523022955079328623947833520$
\newline
$567532298244004442294156837342370982002330874074322076777185746207730$
\newline
$323482406094614280046$
\newline
$+14.23622856322181332287122301085588169299871236208494399568695437825$
\newline
$6069322896396761526007136189745757467102551375667154010366364994538731$
\newline
$7916271113823253751110948972775762217941663830770714262041755664035323$
\newline
$9671078789529204404394764315531582588051352309327292004343654135172820$
\newline
$7780017861238006999109644383198471665302823015355865202971277187847669$
\newline
$1974168218415293165267046606327405458655765280027732495125802150527924$
\newline
$57282410834191507107658393848458313664113623935800293262678700791600125$
\newline
$465010766853i$.
\newline
Let us denote this approximation of $\lambda_1$ as $a$. 
A numerical check  indicates that  
$|\zeta^{\circ 3}(a) -  a|$  agrees with zero to $495$ decimal places.
\newline \newline
Our main reason  for requiring so much precision is that we will be
repeatedly solving equations of the form  $\zeta(u) = v$ for $u$, in each case
replacing  $v$ with the previous $u$, to  construct  lists of
(usually) $100$ elements of a branch of the backward orbit of a nontrivial Riemann
zero, looking for the $u$'s near pre-selected 
fixed points.  As  the procedure 
repeats  $100$ times,  there is an accumulation of numerical error, and
in this situation very high precision  is needed to  maintain enough
accuracy to ``see'' the spirals formed by these branches in our plots.
\section{\sc a tour of A sub phi}
We are interested
in $A_{\phi}$ because 
plots of this set
make visible the underlying
structure of the
network of $\zeta^{\circ n}$ 
pre-images of the critical line for
all $n$ at once:
(1) the nontrivial
 zeros of the
$\zeta^{\circ n}$
lie in  bulbs
of $A_{\phi}$ on filaments $F$ 
decorating
the border of $A_{\phi}$, and 
(2) one
 $\zeta^{\circ n_{_F}}$
pre-image of the critical line
transects each such $F$.
(Claims 1 and 2 are not, of course
logically equivalent;
we are summarizing computer
observations that we will describe
in more detail below.)
Thus the  structure
of union of rooted trees 
visible in plots of $A_{\phi}$
is apparently graph-isomorphic
to a corresponding 
structure for the
point set 
$$
\bigcup_{\Re(z) = \frac 12} 
\zeta^{\circ -}(z) = \mathcal{U}\, \mbox{(say.)}
$$ 
This observation informs our
discussion of the trees $T$ in the 
next section. We pretend 
 that we have stated
 a rigorous
 definition of the decoration notion 
 and  definite conditions
 for the membership of
 a given complex number in a given filament.
 In view of the relationship between
 $\mathcal{U}$
 and  $A_{\phi}$, this
 should not cause problems: each
 filament $F$ may be identified with 
 one (of the many)
  $\zeta^{\circ n_{_F}}$ pre-images of the critical line,
  the definition of which
  could be made precise.
 But we should say explicitly that ``$A$ decorates $B$''
 is a transitive relation and that
  the filaments are subsets of $A_{\phi}$.
\newline \newline
 In Figure 3.1, for example, the 
points at the junctions of four colors 
represent zeros of $\zeta^{\circ 2}$;
the zeros in the long filaments are
nontrivial.
\begin{figure}[!htbp]
\centering
\includegraphics[scale = .4]{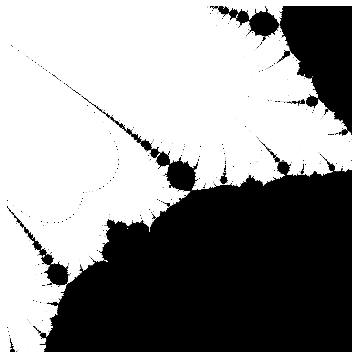}
\hskip .05in
\includegraphics[scale = .4]{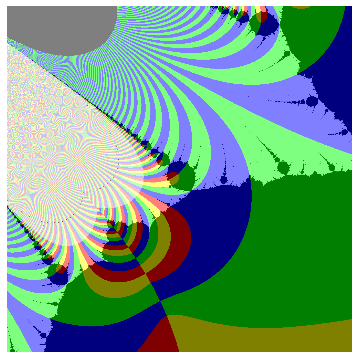}
\vskip .07in
\sc{fig. 3.1: left: $A_{\phi}$ at the edge of the main cardioid; right: superimposed quadrant plot of $\zeta^{\circ 2}$}
\end{figure}
The right panel  of 
Figure 3.2 shows a quadrant
plot of $s \mapsto \zeta(s) - s$
superimposed on a plot of $A_{\phi}$; 
the fixed points of zeta appear as the 
junction of four colors. 
 The left panel
depicts the nontrivial Riemann zeros using
the same scheme (a quadrant plot of zeta.)
\begin{figure}[!htbp]
\centering
\includegraphics[scale=.4]{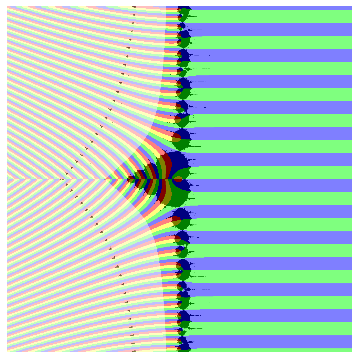}
\hskip .1in
\includegraphics[scale=.4]{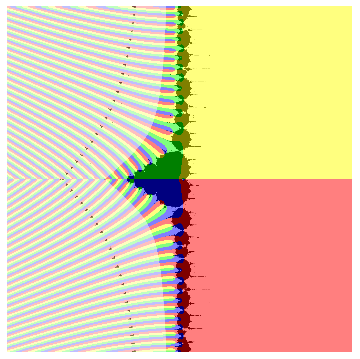}
\vskip .07in
\hskip 0in
\sc{fig. 3.2: left: quadrant plot of zeta; right: quadrant plot of $\zeta(s) -s$; both superimposed on 
$A_{\phi}$}
\end{figure}
\newline \newline 
The filled Julia set of zeta
(the points in $\bf{C}$
with bounded orbit
under iteration by zeta)
is $\bf{C}$ $- A_{\infty}$.
The basin
$A_{\phi}$ appears to be 
dense in  $\bf{C}$ - $A_{\infty}$.
The sets $\bf{C}$ - $A_{\infty}$ 
(Figure 3.3) and $A_{\phi}$, 
regarded as regions in
the complex plane,
 are indistinguishable
in our plots
but they are not identical.
\begin{figure}[!htbp]
\centering
\includegraphics[scale=.85]{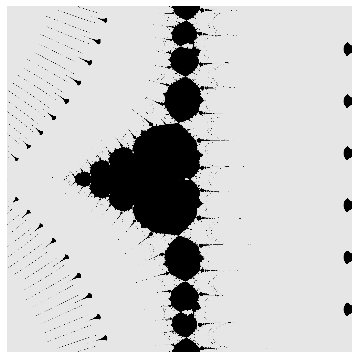}
\vskip .4in
\sc{fig. 3.3:  \bf{C} - $A_{\infty}$}
\end{figure}
For example,
there is an infinite number of real
zeta 
fixed points (\cite{W}, Theorem 1)
that belong to
$(\bf{C}$ $-A_{\infty})-A_{\phi}$.
In addition, there 
appear to be infinite families of 
non-real zeta
$k$-cycles for each integer $k \geq 1$
in $(\bf{C}$ $-A_{\infty})-A_{\phi}$. 
\newline \newline
 Zero lies in $A_{\phi}$ (\cite{W}, Theorem 1.)
This set is a fractal decorated with 
numerous long filaments
(Figure 3.1.)
Zeroes of the $\zeta^{\circ n}$  
lie on the filaments.
Because zero is an element of $A_{\phi}$,
 we know that the whole
backward orbit 
$\zeta^{\circ -}(0)$
lies in $A_{\phi}$.
Because the 
pre-images of  
nontrivial Riemann zeros  under
iterates of zeta lie on the 
filaments,
the itinerary of a point in the backward orbit
of a nontrivial zero $\rho$ visits several 
filaments at the edge
of $A_{\phi}$ before coming to
$\rho$.
(Some but not all pre-images of the trivial
zeros also lie on filaments.)
 \newline \newline
 The set $A_{\phi}$
seems to
comprise \newline 
\newline
(1) a heart-shaped, seven-lobed central body, 
which we will call the main cardioid.
\newline \newline
(2) two major filaments 
of bulbs of various irregular shapes
that emanate from the main cardioid,
transected by the critical line and containing
one nontrivial Riemann zero in each bulb
(right panel of Figure 3.2.)
 \newline \newline
(3) infinitely many  
blunt processes and long filaments 
decorating the main cardioid
and each of the irregular bulbs.
The filaments comprise smaller 
copies of the bulbs,
which, in turn, are decorated with
similar filaments,
\it ad infinitum. \rm 
(Figure 3.1.) 
Thus, when we plot them, the set of 
filaments decorating $A_{\phi}$ exhibit
a visible tree structure.
\newline \newline
The visible features 
described in (1) - (3) were evident in
Woon's plots of $\bf{C}$ $- A_{\infty}$ 
\cite{W}.
The filaments appear to be 
zeta-iterate pre-images (close copies) of
the two major filaments. 
Pre-images of the real axis pass through
the blunt processes and contain
pre-images of the trivial
zeros. For example, 
the right panel of Figure 3.1,
superimposes a quadrant plot 
of $\zeta^{\circ 2}$ on the
left panel, so that
junctions 
of four differently-colored regions
each represent a zero of $\zeta^{\circ 2}$. There are
three long filaments depicted in this image
containing zeros, the immediate zeta images 
of which are 
nontrivial Riemann zeros; but between the lower two
such filaments is a blunt
process transected by a zeta pre-image of the real axis,
and we can see another series of 
$\zeta^{\circ 2}$-zeros lying along this curve.
These are zeta pre-images of the negative
even numbers.
\newline \newline
(4) at each trivial zero $< -18$,
a microscopic, more or less distorted 
copy (zeta-iterate pre-image)
of the entire assemblage 
described in (1) - (3).
(By ``microscopic'' features
we mean features so small
that they
can only be visualized by a 
change of scale from that of Figure 3.3.)
In Figure 3.4, we show copies
of the main cardioid
near the trivial zeros $-28, -26, -24, -22$
 superimposed on quadrant plots of 
$\zeta^{\circ 2}$ in the same squares.
\begin{figure}[!htbp]
\centering
\includegraphics[scale=.45]{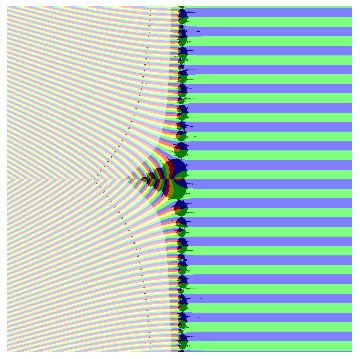}
\hskip .05in
\includegraphics[scale=.45]{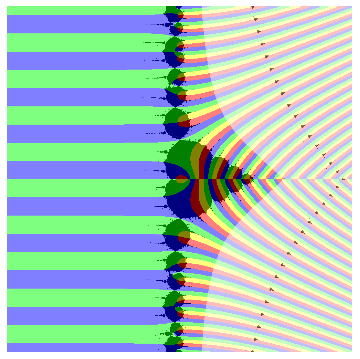} 
\vskip .05in
\includegraphics[scale=.45]{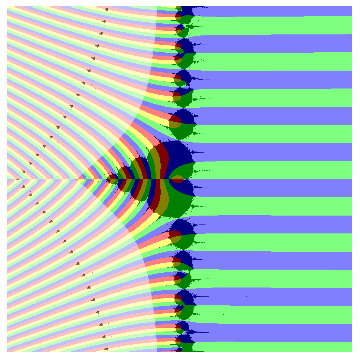} 
\hskip .05in
\includegraphics[scale=.45]{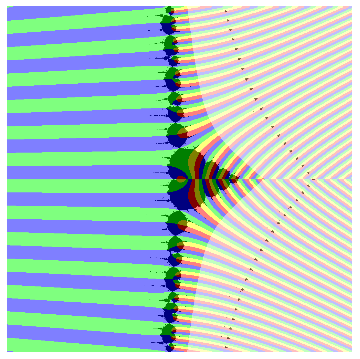} 
\vskip .1in
\sc{fig. 3.4:  $A_{\phi}$ copies near
 $s = -28, -26, -24, -22$  superimposed on $\zeta^{\circ 2}$ 
quadrant plots}
\end{figure} 
The size of these features decays exponentially
with distance from zero. Their left-right orientation
alternates.
 We speculate that the alternation
can be derived from the 
alternating sign of the real
derivative $\frac {d \, \zeta(x)}{dx}|_{x = -2n}$. 
\newline \newline
Because these copies exist
on the left half of the real axis,
its zeta pre-images also contain
complete copies of $A_{\phi}$. The upper left panel of
Figure 3.5 depicts the first bulb in the major filament
in the upper half plane. It is a $10$ by $10$
square centered at $\rho_1$. Along its
border we see 
an apparently infinite set of filaments alternating with
an apparently infinite set of blunt processes.
\begin{figure}[!htbp]
\centering
\vskip 0.1in 
\hskip .2in
\includegraphics[scale=.45]{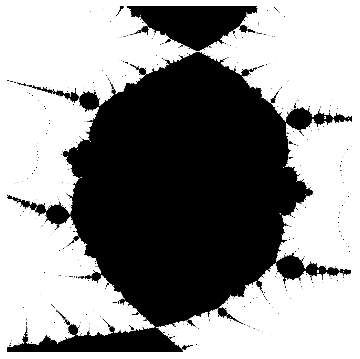} 
\includegraphics[scale=.45]{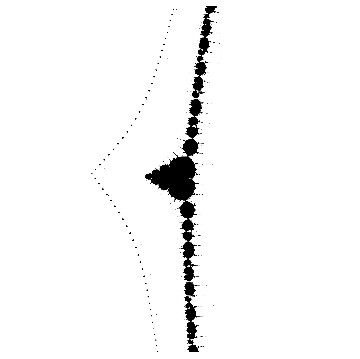} 
\vskip .1in
\hskip .2in
\includegraphics[scale=.45]{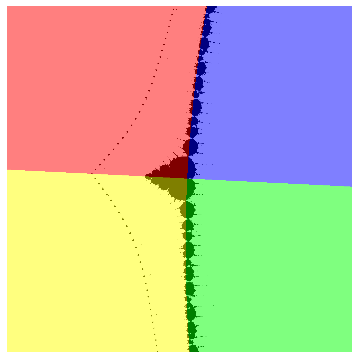}
\includegraphics[scale=.45]{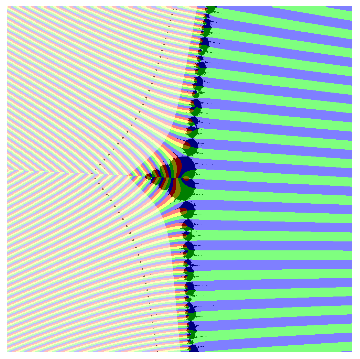}
\vskip .2in
\hskip 1.1in 
\sc{fig. 3.5: a copy of $A_{\phi}$ near a blunt process of $A_{\phi}$; left: with superimposed quadrant plot of $\zeta^{\circ 3}$; right: superimposed quadrant plot of $\zeta^{\circ 4}$}
\end{figure} 
\begin{figure}[!htbp]
\centering
\includegraphics[scale=.45]{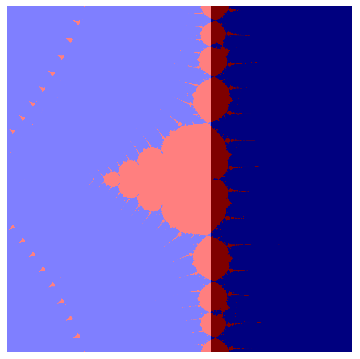} 
\includegraphics[scale=.45]{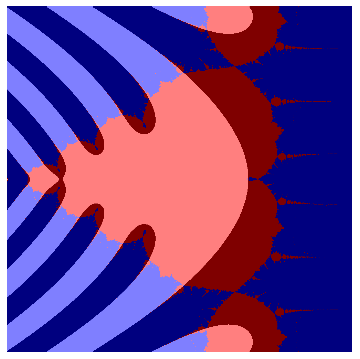} 
\vskip .2in
\includegraphics[scale=.45]{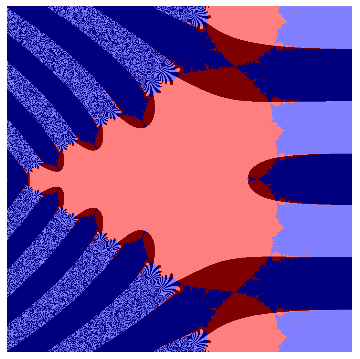}
\includegraphics[scale=.45]{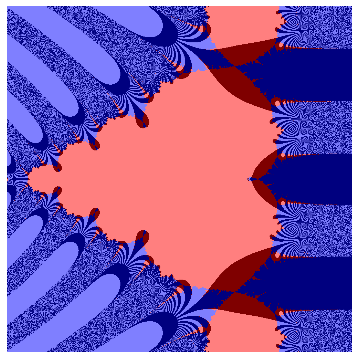}
\newline
\hskip .2in 
\sc{fig. 3.6: darker regions are $\zeta^{\circ n}$ 
pre-images of $\Re(s) > \frac 12$ near the main
cardioid; $n = 1$ in row 1 column 1; $n=2$ in row 1 column 2; $n=3$ in row 2 column 1;  $n = 4$ in row 2 column 2}
\end{figure} 
Our tests demonstrate that the filaments are
transected by $\zeta^{\circ n}$ 
pre-images of the critical line 
for $(n = 1, 2, 3, ...)$,
and that the blunt processes are 
transected similarly by
$\zeta^{\circ n}$ pre-images of the real axis. 
The other three panels depict a small 
copy of $A_{\phi}$ to the right of
the largest
blunt process on the right side of the bulb 
shown in the upper left panel.
In the lower panels, a quadrant plot of $\zeta^{\circ 3}$ 
in the left panel and of $\zeta^{\circ 4}$ 
in the right panel have been 
superimposed upon this copy.
Evidently, it is
a $\zeta^{\circ 3}$ pre-image of $A_{\phi}$.
\newline \newline
(5) Our observations indicate that for
 each filament $F$ decorating $A_{\phi}$
 there is a positive integer $k_F$
 (say, the degree of $F$)
 such that each bulb of $F$ contains 
 one 
 nontrivial $\zeta^{\circ k_F}$ zero,
 and no nontrivial zeros of  $\zeta^{\circ k}$
 for any $k \neq k_F$. 
 Even under the Riemann hypothesis, it would not be
necessary from
 first principles that degree $k$ filaments
 are transected by $\zeta^{\circ k -1}$ 
 pre-images of the critical line,
even though that is the simplest possibility. 
But it seems to be
 the case. In Figure 3.6, the $\zeta^{\circ k -1}$
 pre-images of the critical line transecting 
 degree $k$ filaments decorating
the main cardioid are shown for 
$k = 1, 2, 3$ and $4$.
\section{\sc branches interpolated by spirals}
\begin{figure}[!htbp]
\centering
\includegraphics[scale=.33]{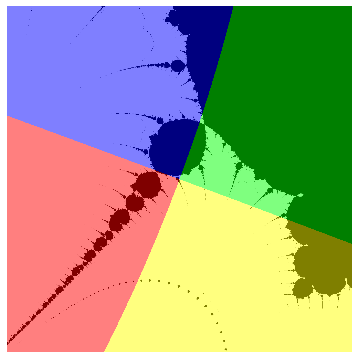} 
\hskip .7in
\includegraphics[scale=.33]{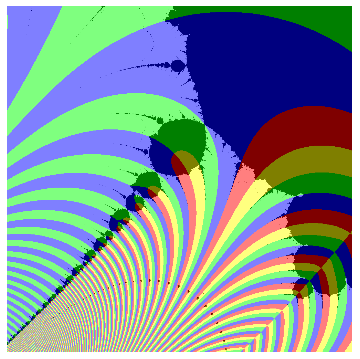} 
\vskip .1in
\includegraphics[scale=.33]{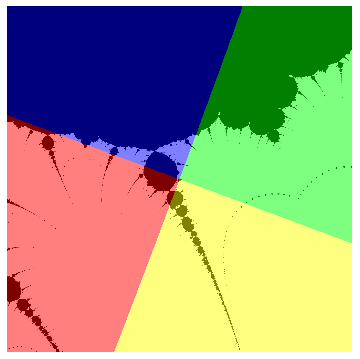}
\hskip .7in
\includegraphics[scale=.33]{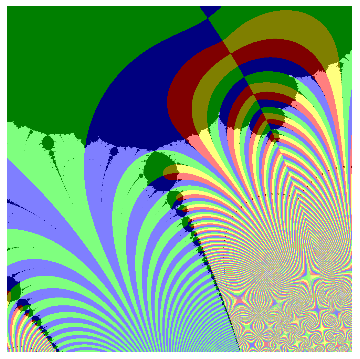}
\vskip .1in
\includegraphics[scale=.33]{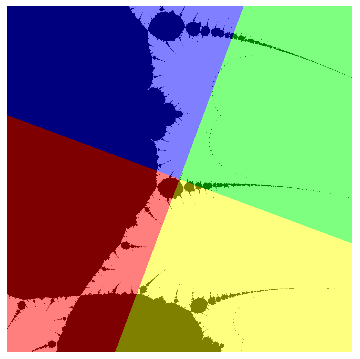}
\hskip .7in
\includegraphics[scale=.33]{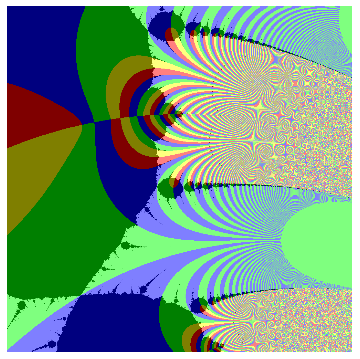}
\vskip .1in
\hskip .23in
\includegraphics[scale=.33]{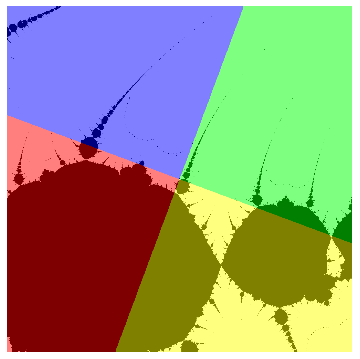}
\hskip .7in
\includegraphics[scale=.33]{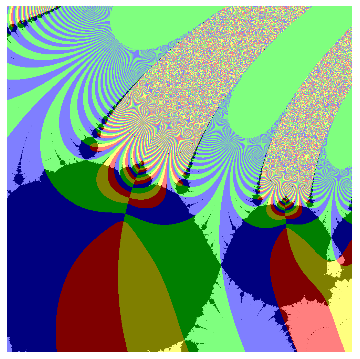}
\newline 
\vskip .1in
\sc{fig. 4.1: zooming in on $\psi_{\rho_1}$; column 1: quadrant plots of $\zeta(s) - s$; column 2 (top to bottom): quadrant plots of $\zeta^{\circ n}, n = 3, 4, 5, 6$}
\end{figure} 
For any integer $L > 0$ 
there appears
to be an infinite set of
zeta cycles $\Lambda = ( \lambda_0, ..., \lambda_{L-1})$
that pick out a linearly ordered subset
 $B_{\rho, \Lambda} = (a_0, a_1, a_2, ...)$
 of $\zeta^{\circ -}(\rho)$ 
such that (1) $a_0 = \rho$,
(2) $a_n = \zeta (a_{n+1}), n = 0, 1, 2, ....$
and (3) for each $j = 0, 1, 2, ..., L-1$,
the subsequences
 $b_j = (a_j, a_{j + L}, a_{j + 2L}, ...)$
appear to converge to  $\lambda_j$.
The sequence $b_j$ is a branch of $\zeta_L^{\circ -}(a_j)$.
In most of the cases 
that we have examined,
each $b_j$ appears to be
interpolated by a
spiral $s_{\rho_N, \, \lambda_j}$
with center $\lambda_j = \lim b_j$. 
The $\lambda_j$ are repelling
 fixed points of $\zeta_L$. 
\newline \newline
 Now we offer (speaking loosely) a geometric description
 of some of the  $b_j$ in terms of the
 basin of attraction $A_{\phi}$. 
 (It applies to most, but not
 all, instances we have examined to date.)
 A variety of filaments decorate 
 $A_{\phi}$, but here we restrict attention to those 
 that decorate the main
 cardioid. We assign the set of filaments 
 a structure of union of 
rooted trees $T$ as follows.
A filament $F \in T$
 is the parent of a filament $G \in T$
 if and only if $G$ decorates $F$ and there is
 no intermediate filament $H \in T$ such that
 $G$ decorates $H$ and $H$ decorates $F$.
 The filaments containing
 nontrivial Riemann zeros have no ancestors, 
 but they are not unique in this respect.
  \newline \newline
Now fix integers
$L \geq 1, N \neq 0$. There is an infinite
set of zeta $L$-cycles 
$\Lambda = (\lambda_0, \lambda_2, ..., \lambda_{L-1})$
such that for each integer 
$\Delta = 0, 1, 2, ..., L-1$,
there is  a tree 
$T_{\Delta}$ of filaments decorating the main 
cardioid of $A_{\phi}$, and a 
path $P_{\Delta}  = (F_0, F_1, ...)$ in $T_{\Delta}$
 with 
 $k_{F_m} = \Delta + m L$
 and such that
 (if $m > 0) \, F_m$
 decorates the $|N|^{\rm th}$ bulb of
 its parent $F_{m-1}$. 
 As in the first column
 of Figure 4.1,
 the filaments  in
 $P_{\Delta}$ spiral around  $\lambda_{\Delta}$.
In our graphic visualizations,
the apparent size of $F_m$
  decays 
 exponentially with $m$.
 Something like this would seem to be a
 necessary condition of the relation
 $\lambda_j = \lim b_j$
 we mentioned above.
 \newline \newline
Each bulb of $F_m$ contains
a nontrivial zero $w$ 
of $\zeta^{\circ \Delta + mL}$ and
$\zeta^{\circ \Delta + mL-1}(w)$ is
a nontrivial Riemann zero. Which one?
Let $w_{m, N}$ be
the nontrivial $\zeta^{\circ \Delta + mL}$ zero
belonging to the $|N|^{\rm th}$
 bulb of
 the filament $F_m$  in $P_{\Delta}$. 
  For $m > 0$,
 $\zeta_L(w_{m, N}) = w_{m - 1, N}$,
 so the sequence 
 $(w_{0, N}, w_{1, N}, ...)$
 is a branch of $\zeta_L^{\circ -}(w_{0, N})$. 
 In our observations, the 
 zeta image of bulb $|N|$ of a filament
 $F$ in $T_{\Delta}$
 with $k_F > 1$
 is bulb $|N|$ of its parent filament
 in $T_{\Delta}$.
Therefore
 $\zeta^{\circ \Delta + mL -1}(w_{m, N}) = \rho_{\pm N}$. 
 \newline \newline
 The observation that
  $\zeta_L(w_{m, N}) = w_{m - 1, N}$
   suggests that
 there should be graph isomorphisms between
 subgraphs of the rooted tree graphs
 associated to the $\zeta_L^{\circ -}$ on
 one side and subgraphs of the  
 trees $T$ decorating the main cardioid of 
 $A_{\phi}$ on the other.
 We should mention that the situation 
 for \it copies \rm of the main 
 cardioid such as the ones
 illustrated in Figures 3.4 and 3.5
 is different; except to say that
 the zeta images of copies are also copies,
 we will not discuss it
 further in the present article.
\section{\sc spirals interpolating a branch of 
the backward zeta orbit of a Riemann zero}
\subsection{Single spirals interpolating a branch.}
When $L = 1$, $\Lambda  = \{\lambda_1 \}$ where
$\lambda_1$ is a repelling
zeta fixed point; there appear to be 
at least three categories of such  points:
$\psi_{-2n}$ (say) lying near the 
trivial zeros $-2n = -20, -22, ...$; 
zeta fixed points $\psi_{\rho^*}$  near each nontrivial
Riemann zero $\rho^*$, 
and eight
fixed points lying at the boundary of 
the main cardioid (right panel, Figure 3.2.)
How near? In the case of 
the $\rho^*$, one can form an impression by
keeping in mind that this figure
depicts a $120$  by $120$ square (section 8.) 
The distances $|-2n -\psi_{-2n}|$ are a great dealer smaller;
we omit the details. 
\newline \newline
There are exactly two filaments $F$ 
with degree $k_F = 1$ decorating the main cardioid; 
one of them  contains
$\rho_N$  in its  $|N|^{\rm th}$ bulb
$ = \beta_N$, say. 
Our observations are consistent
with the following proposition.
The point $\psi_{\rho_{_N}}$
lies at the border of 
the $|N|^{\rm th}$
bulb of
a filament $F^*$ 
with $k_{F^*} = 2$ decorating
$ \beta_N$.
This ramifies: if $\psi_{\rho_N}$
 lies at the border of 
the $|N|^{\rm th}$
bulb of
a filament $F^*$ 
then there is a child filament
$F'$ of $F^*$ such that
 $\psi_{\rho_N}$
 lies at the border of 
the $|N|^{\rm th}$
bulb of $F'$ 
and  $k_{F'} = k_{F^*}+1$.
\newline \newline
This is illustrated by  the left column of
Figure 5.1 for $N = 1, 2, 3, 4$. It depicts
quadrant plots
of $s \mapsto \zeta(s) - s$, so that $\psi_{_{\rho_N}}$
shows up as four-color junctions on 
the depicted filament (say, $F_2$.)
The right column shows
 quadrant 
 plots of 
 $\zeta^{\circ 3}, \zeta^{\circ 4}, \zeta^{\circ 5}, \zeta^{\circ 6}$
 in rows 1, 2, 3 and 4, respectively,
 all superposed on plots of $A_{\phi}$.
  The squares have side length
 $.2, .02, .002,$ and $.0002$ 
 in rows 1, 2, 3 and 4, respectively.
 The center of the squares in row $N$ is 
 $\psi_N$, so the panels are 
 depicting the region around this
 point at smaller and smaller scales.
\begin{figure}[!htbp]
\centering
\includegraphics[scale=.33]{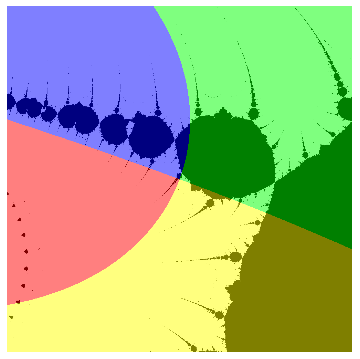} 
\hskip .7in
\includegraphics[scale=.33]{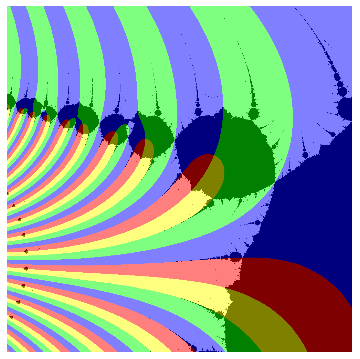} 
\vskip .1in
\includegraphics[scale=.33]{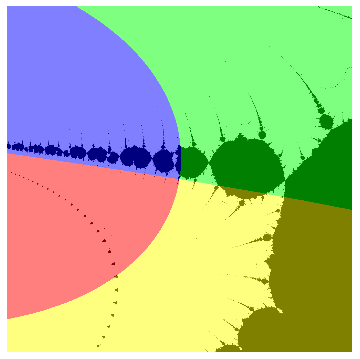}
\hskip .7in
\includegraphics[scale=.33]{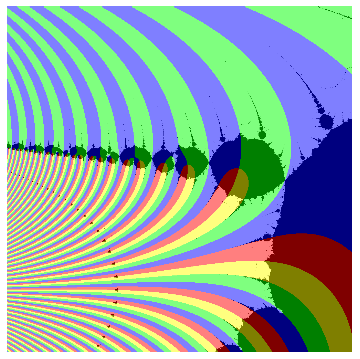}
\vskip .1in
\includegraphics[scale=.33]{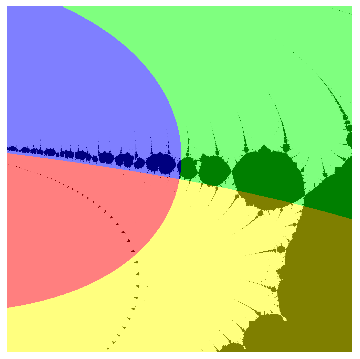}
\hskip .7in
\includegraphics[scale=.33]{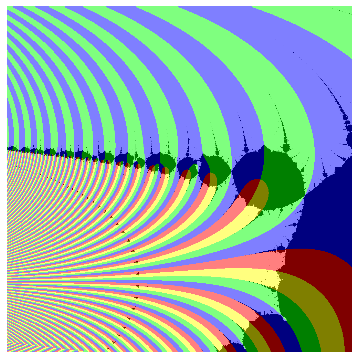}
\vskip .1in
\hskip .23in
\includegraphics[scale=.33]{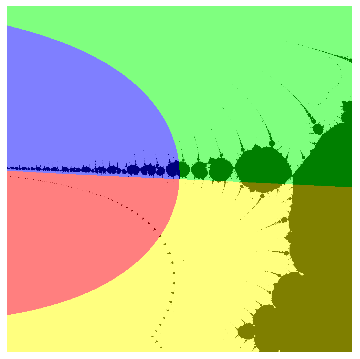}
\hskip .7in
\includegraphics[scale=.33]{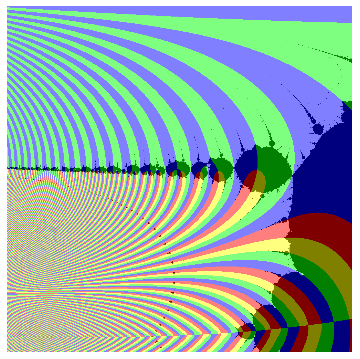}
\newline 
\vskip .1in
\sc{fig. 5.1: top to bottom: $\psi_{\rho_n}, 1 \leq n \leq 4$; column 1: quadrant plots of $\zeta(s) - s$; column 2: quadrant plots of $\zeta^{\circ 2}$}
\end{figure}
\newline \newline
Figure 4.1 displays the original indications
we had that some branches of
the inverse of zeta
lie on spirals. 
It zooms in on the
illustration of $\psi_{\rho_1}$
in the left panel of 
the top row of Figure 5.1. 
The center of that panel is the fixed point,
lying on the border
of a filament $F_3$ (say) 
decorating the lower border of 
the largest full bulb.
In the right 
panel of the top row of Figure 4.1,
a quadrant plot of 
$\zeta^{\circ 3}$ has been superimposed
on a plot of the corresponding region of
$A_{\phi}$; we see that $F_3$ contains zeros
of $\zeta^{\circ 3}$. (Simple tests 
show that they are
nontrivial zeros in the sense of
the introduction.) 
In the lower rows,
the zoom is repeated and $\psi_{\rho_1}$
is seen to lie near a
still-smaller filament decorating the
bulb near the center of the
figure just above it. The right column
depicts quadrant plots of $\zeta^{\circ 4}$,
$\zeta^{\circ 5}$ and $\zeta^{\circ 6}$
for the squares opposite them
in the left column. So in Figure 4.1 we are
seeing zeros of these functions (again
nontrivial.)
They also appear 
(in virtue of the shapes of the
underlying $A_{\phi}$-bulbs) 
to be
$\zeta^{\circ 3},
\zeta^{\circ 4}, \zeta^{\circ 5}$
pre-images of nontrivial Riemann zeros.
The rapid reduction of scale from one row
to the next attests to a similar
reduction of the distances of these 
pre-images from $\psi_{\rho_1}$
(which, as we have remarked, is not surprising.)
The possibility that they may be traveling
on spirals emerges from a look at the angles
that the filaments $F_3$ and (say) $F_4, F_5, 
F_6$ make with the horizontal. 
These observations led us to do the numerical
tests described in the last section.
\newline \newline
We made a survey of the spirals $s_{\rho,\psi_{\rho^*}}$ 
for various choices of $\rho$ and $\rho^*$. 
We  made a table of $\psi_{\rho_n}, 1 \leq n \leq 100$ 
with $500$ digits of precision
and we used a table of nontrivial Riemann zeros 
with $300$ digits of precision
made by Andrew Odlyzko \cite{O}. We made tables of
the $z_k$ in various 
$B_{\rho, \psi_{\rho^*}}$ to high precision,
proceeding inductively. 
We set $z_1 = \rho$ and, given 
a value of $z_k$, after using
\it Mathematica's \rm FindRoot command 
to solve $\zeta(s) = z_k$ in the vicinity of $\psi_{\rho^*}$,
we set $z_{k+1}$ equal to the solution. We began with $300$ 
$z_k$ for each $B$ and used tests of reliability of
each $z_k$ to truncate the list; typically, 
we ended up with a least $100$
consecutive $z_k$. 
\newline \newline
We use polar coordinates $(r(z), \theta(z))$
to denote a typical point $z$ on $s_{\rho,\psi_{\rho^*}}$
such that (1) $r(z) = |z - \psi_{\rho^*} |$,
(2)  $\theta(z)$ is chosen so that 
$\theta(z)  
\equiv \arg(z - \psi_{\rho^*}) \pmod {2\pi}$,
and (3) $\theta(z)$ varies continuously and monotonically 
as $z$ moves around the spiral in a fixed direction. 
In other words, $\theta(z)$ 
behaves up to a  mutiplicative 
constant like a winding number.
Then $r(z)$ appears to decay exponentially with $\theta(z)$.
(Of course we are only able to check this
for $z \in B_{\rho,\psi_{\rho^*}}$, that is,
for the $z_k$ we propose are interpolated by 
$s_{\rho,\psi_{\rho^*}}$, because no other test for 
membership in $s_{\rho,\psi_{\rho^*}}$ is  available
to us.)
Therefore it was not practical to plot the spirals 
$s_{\rho,\psi_{\rho^*}}$ without re-scaling, 
so we plotted the points 
$(\log r(z), \theta(z))$ instead.
This procedure 
 everts the apparent spirals: 
if $k, j$ are such that  
$r(z_k) < 1$ and $r(z_j) < 1$, then 
$r(z_j) <  r(z_k)$ implies
that the plotted point 
$(\log r(z_j), \theta(z_j))$
is further from
 the center of the re-scaled interpolating spiral
than the point 
$(\log r(z_k), \theta(z_k))$: the
reverse of the situation before re-scaling. 
(The direction of winding
 of the spiral is also reversed because the 
 logarithms take negative values.) The
   points near the center of the re-scaled spiral
   depict $z_k$ for smaller values of $k$ 
   for which  $z_k$ is closer to
   $\rho$ and further from $\rho^*$.
  They are crowded so closely, in spite of our
   re-scaling, that
   the interpolating curve near $\rho$ is 
   obscured.\newline \newline 
   \begin{figure}[!htbp]
\vskip 0.1in 
\hskip .7in
\includegraphics[scale=.3]{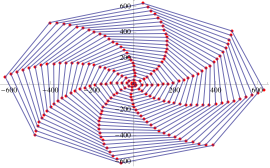} 
\hskip .7in
\includegraphics[scale=.3]{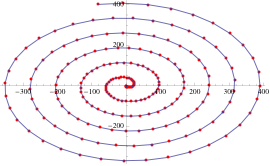} 
\vskip .1in
\hskip .7in
\includegraphics[scale=.3]{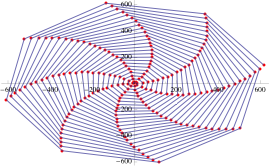}
\hskip .7in
\includegraphics[scale=.3]{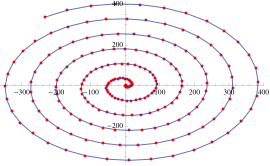}
\vskip .1in
\hskip .7in
\includegraphics[scale=.3]{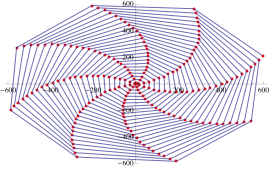}
\hskip .7in
\includegraphics[scale=.3]{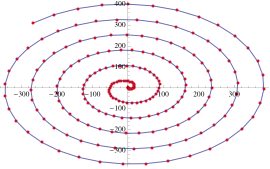}
\vskip .1in
\hskip .7in
\includegraphics[scale=.3]{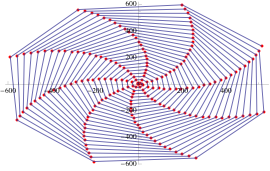}
\hskip .7in
\includegraphics[scale=.3]{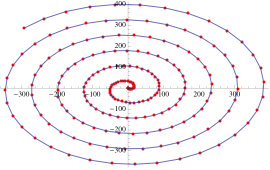}
\vskip .1in
\hskip .7in
\includegraphics[scale=.3]{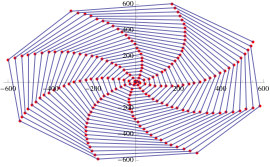}
\hskip .7in
\includegraphics[scale=.3]{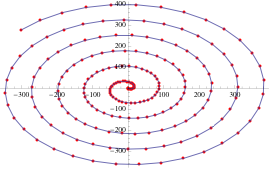}
\newline 
\vskip .1in
\hskip 0in
\sc{fig. 5.2: backward orbit branches for 
$\rho_n \, (1 \leq n \leq 5)$ 
centered on two cardioid zeta fixed points}
\end{figure}
Figure 5.2 depicts branches of backward orbits of 
   $\rho_n (1  \leq n  \leq 5)$
   spiraling around two fixed point 
   $\approx -14.613 + 3.101i$
   (left column) and  
   $\approx  - 5.28 + 8.803i$
   (right column)
   on the
   border of the main cardioid; 
   we omit the 500 digit decimal
   expansions, which are easy to compute 
   using \it Mathematica\rm's 
   FindRoot command.
   Figure 5.3 depicts branches of
   backward orbits of $\rho_n$ 
   spiraling around $\psi_{\rho_{_n}} (1 \leq n \leq 10)$.
   (We omit their precise expansions  for the same reason.)
\begin{figure}[!htbp]
\centering

\includegraphics[scale=.3]{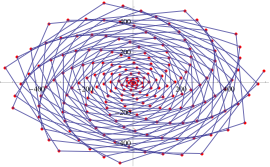} 
\hskip .7in
\includegraphics[scale=.3]{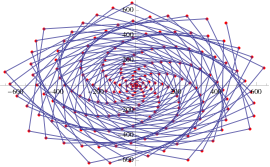} 
\vskip .1in

\includegraphics[scale=.3]{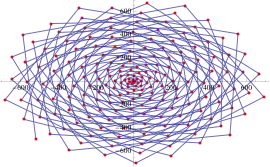}
\hskip .7in
\includegraphics[scale=.3]{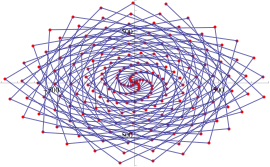}
\vskip .1in

\includegraphics[scale=.3]{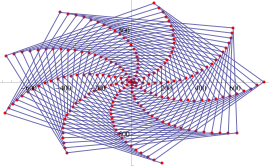}
\hskip .7in
\includegraphics[scale=.3]{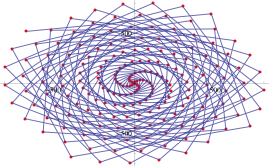}
\vskip .1in

\includegraphics[scale=.3]{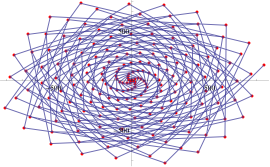}
\hskip .7in
\includegraphics[scale=.3]{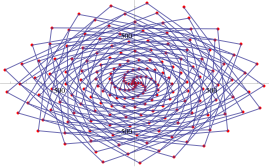}
\vskip .1in
\hskip .4in
\includegraphics[scale=.3]{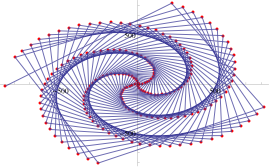}
\hskip .7in
\includegraphics[scale=.3]{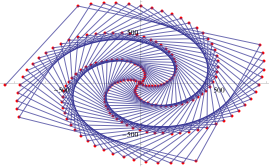}
\newline 
\vskip .1in

\sc{fig. 5.3: backward orbit branches 
for $\rho_n$ near $\psi_{\rho_n}, 1
\leq n \leq 10$}
\end{figure}
\subsection{An example.}
\begin{figure}[!htbp]
\centering
\includegraphics[scale=.4]{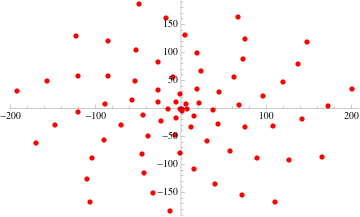}
\hskip .1in
\includegraphics[scale=.4]{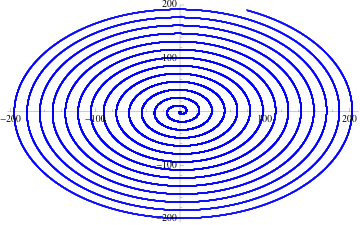}
\vskip .07in
\includegraphics[scale=.4]{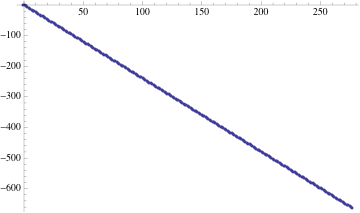}
\hskip .1in
\includegraphics[scale=.4]{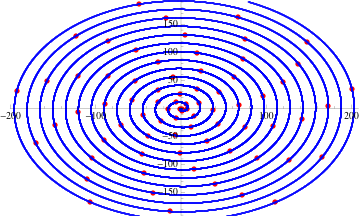}
\vskip .07in
\sc{fig. 5.4: logarithmic scaling of  
$s_{\rho_1,\hskip .03in \rho_1}$}
\end{figure}
The first panel of Figure 5.4 plots the point set
$B_{\rho_1,\psi_{\rho_1}}$
(re-scaled as described, 
and shifted to place the apparent spiral's 
center at the origin.)
It offers the 
appearance that the $z_k$ (in red) form arms something
like those of
a spiral galaxy;
this seems to be a result of
nearly regular
growth of $\theta(z_k)$ with $k$. 
But $z_k$ for consecutive $k$ do
not lie in adjacent positions on these arms;
in the second panel,  
the $z_k$ are connected by chords in the same
order as they appear in the sequence 
$B_{\rho_1, \psi_{\rho_1}}$: 
vertices $v, w$  representing $z_v, w = \zeta(z_v)$
are connected by a chord.
The lower left panel of Figure 5.4 is a plot 
of $\log |z_k - \psi_{\rho_1}|$ vs. $k$;
it is clear that the distances of the $z_k$
(colored red)
from the fixed point $\psi_{\rho_1}$ at the center
of the spiral are decaying exponentially.
The lower right panel of Figure 5.4 
depicts a spiral curve (colored blue)
that approximately interpolates the
 $z_k$; we found this curve using
the NonLinearModelFit command in 
\it Mathematica\rm.
The equation of the curve is 
$$ \log r(z) = a + b \theta(z) +  c \exp(d \theta(z)) ,
$$
with 
$$a \approx 0.05575203301551956560399459579161529353,$$ 
$$b \approx -2.39481894384498740085074310912697832305,$$
$$c \approx -2.8680355917721941635331399485184884 \times 10^{-120},$$
and
$$d \approx 0.97375124237020440301256901292961731822.$$
The absolute value of $c$ is so small that 
this is quite close to being
the equation of a logarithmic spiral.
In the section on error terms below
we will compare directly the loci of the $z_k$ with 
logarithmic spirals.
\newline \newline
In Figure 5.5, we study the 
variation in $\theta(z_k)$. The
left panel plots $\delta_k = \theta(z_{k+1}) - \theta(z_k)$
against $k$ and shows that the $\theta(z_k)$ are
very nearly periodic in $k$. 
\begin{figure}[!htbp]
\centering
\includegraphics[scale=.4]{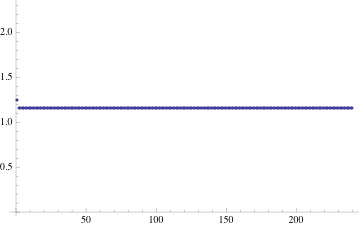}
\hskip .01in
\includegraphics[scale=.4]{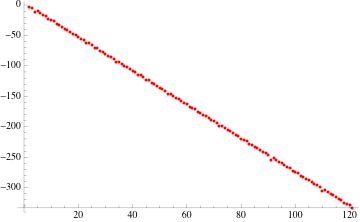}
\vskip .07in
\sc{fig. 5.5: left:  $\delta_k$ vs. $k$ 
right: $\log | \delta_{k+1} - \delta_k|$ 
vs. $k$.}
\end{figure}
The right panel,
which plots 
$\log |\delta_{k+1} - \delta_k|$ against $k$, 
shows that the departure from periodicity
in $\theta(z_k)$ actually appears to 
decay exponentially
with $k$.
However for other choices of 
$\rho$ and $\rho^*$
this no longer holds, and so it is an 
open question whether or 
not it would hold even in this example 
for very large $k$, that is, very close 
to the center of the spiral. 
We remark that the $z_k$ could, of course,
be distributed along a nearly-logarithmic spiral
while also being distributed in a completely
irregular or at least non-periodic way in the
theta aspect, so the two questions are at least
superficially independent.
\newline \newline
Now suppose $(r_1, \theta_1)$ and $(r_2, \theta_2)$ 
lie on a true logarithmic spiral 
$\log r = a + b\theta$. The constants $a$, $b$ are 
determined by any two points of the spiral, hence,
if two pairs of points determine
different values of $a$ and $b$,
then the  curve that the three (or four)
points comprising the pairs lie on 
is not a logarithmic spiral.
We used this idea to test 
$B_{\rho_1,\psi_{\rho_1}} = \{z_1, z_2, z_3, ... \}$  
for the property of being 
interpolated by a logarithmic spiral. 
We performed the
test by solving for $a$ and $b$ 
using the pairs  $(z_1, z_k), k = 2, 3, ....$
\begin{figure}[!htbp]
\centering
\includegraphics[scale=.3]{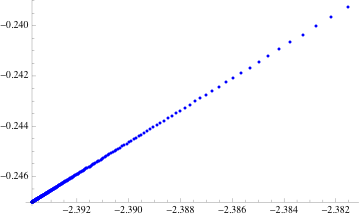}
\hskip .01in
\includegraphics[scale=.3]{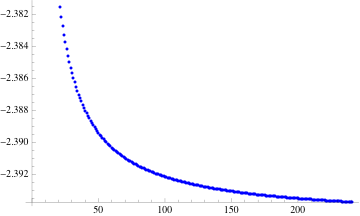}
\includegraphics[scale=.3]{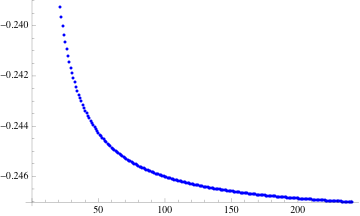}
\vskip .07in
\sc{fig. 5.6: parameters $a, b$ for logarithmic 
curve $\log r = a + b \theta$ induced
 by successive 
points of $B_{\rho_1,\rho_1}$. 
left: $a$ vs. $b$; 
center: $b$ vs. $k$; right: $a$ vs. $k$}
\end{figure}
In Figure 5.6 we have plotted the resulting 
values of $a$ against $b$,
$b$ against $k$, and $a$ against $k$.
Evidently $a$ is roughly linear in $b$,
 and both $a$ and $b$
appear to converge as $k$ grows
without bound. Thus the interpolating 
curve is not a
logarithmic spiral, for then $a$ and $b$ would 
be constants.
But the convergence of $a$ and $b$ 
suggests that the 
interpolating spiral $s_{\rho_1,\psi_{\rho_1}}$ 
resembles a logarithmic spiral 
more and more closely
as it winds inward towards $\psi_{\rho_1}$.
\subsection{Backward orbits near the trivial zeros.}
There appear to be
real zeta fixed points $\psi_{-2n}$
near each trivial zero $-2n \leq -20$.
Whether they lie slightly to the right or 
to the left
of $-2n$ along the real axis appears to depend
upon the parity of $n$. This reflects
 the alternating left-right orientations
of copies (zeta pre-images) of the 
basin of attraction $A_{\phi}$ we see
in Figure 3.4. 
\newline \newline
A branch $B_{\rho,\psi_{_{-2n}}}$ 
of the backward orbit of each nontrivial
Riemann zero $\rho$ lies on a 
curve appearing to pass through or terminate at 
 $\psi_{-2n}$: if $-2n \equiv 0$ (mod $4$),
 then the curve appears to terminate at 
 $\psi_{-2n}$. These curves closely resemble
 straight line segments.  Error terms
 are discussed in section 6. If  $-2n \equiv 2$ (mod $4$), 
 then  (supposing, for the moment, that the curve
 really is a line segment)
 $\psi_{-2n}$ lies near its midpoint.
\begin{figure}[!htbp]
\centering
\vskip 0.1in 
\includegraphics[scale=.4]{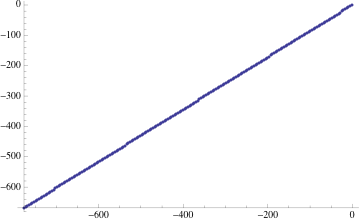} 
\hskip .1in
\includegraphics[scale=.4]{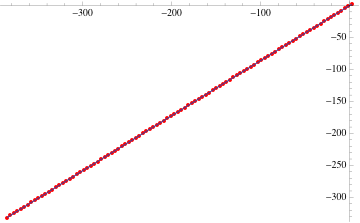} 
\vskip .1in
\includegraphics[scale=.4]{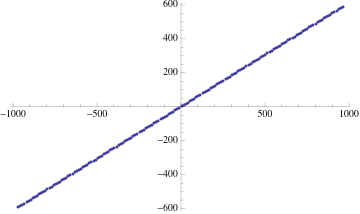}
\hskip .1in
\includegraphics[scale=.4]{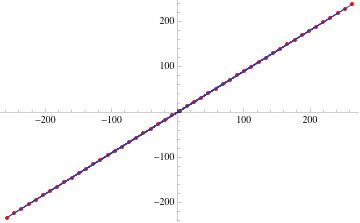}
\vskip .1in
\includegraphics[scale=.4]{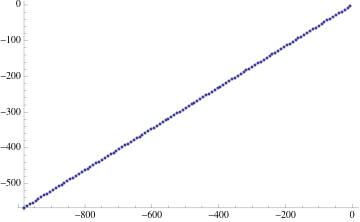} 
\hskip .1in
\includegraphics[scale=.4]{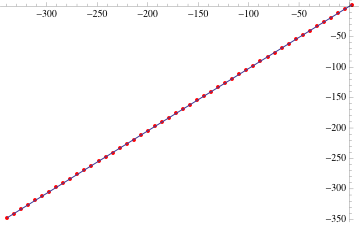} 
\vskip .1in
\includegraphics[scale=.4]{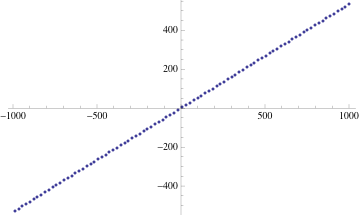}
\hskip .1in
\includegraphics[scale=.4]{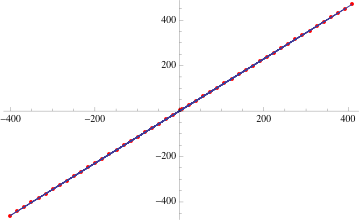}
\vskip .1in
\hskip 0in
\sc{fig. 5.7: column 1: backward orbits of $\rho_1$ near 
$\psi_{-2n}, 10 \leq n \leq 13$}
\sc{column 2: backward orbits of $\rho_{n-9}$ near 
$\psi_{-2n}, 10 \leq n \leq 13$}
\end{figure}
 In Figure 5.7, 
  several of the backward orbits are depicted, 
  re-scaled logarithmically as above.
\newline \newline
This observation is  consistent with the
hypothesis that $B_{\rho,\psi_{-2n}}$ is
interpolated by a spiral such that the
$a_k \in B_{\rho,\psi_{-2n}}$
satisfy $|\arg(a_k-\psi_{-2n}) - \arg(a_{k + 1}-\psi_{-2n})| \approx 2\pi$
for all $k$
if $-2n \equiv 0$ (mod $4$), or 
$|\arg(a_k-\psi_{-2n}) - \arg(a_{k + 1}-\psi_{-2n})| \approx \pi$
if $-2n \equiv 2$ (mod $4$). We discuss this
further in sections 5.5 and 6.2.
\subsection{Several spirals that 
together interpolate a branch.}
\begin{figure}[!htbp]
\centering
\vskip 0.1in 
\hskip 0in
\includegraphics[scale=.4]{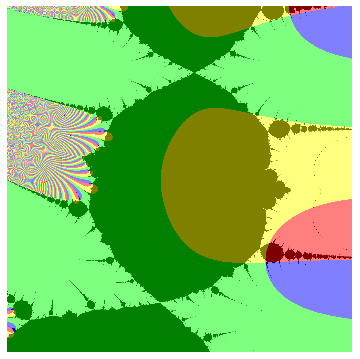} 
\includegraphics[scale=.4]{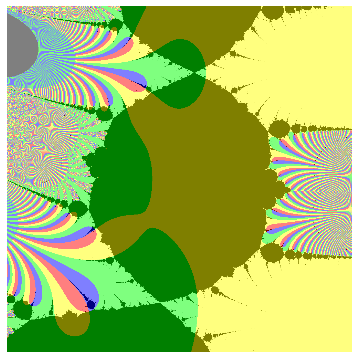} 
\vskip .05in
\includegraphics[scale=.4]{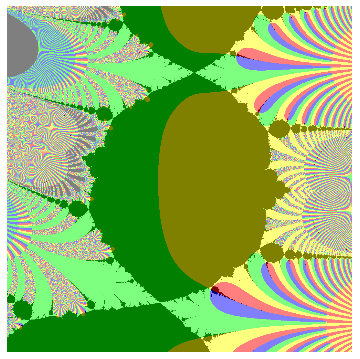}
\includegraphics[scale=.4]{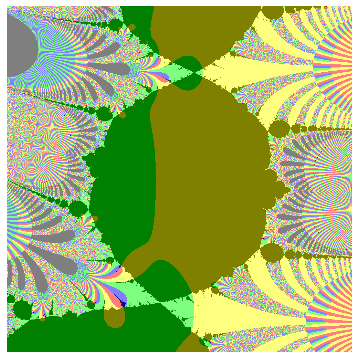}
\hskip 0in 
\vskip.1in
\sc{fig. 5.8: quadrant plots of 
$s \mapsto \zeta^{\circ n}(s) - s, 2 \leq n \leq 5$, near the bulb around $\rho_1$}
\end{figure} 
Figure 5.8 depicts members of zeta $n$-cycles 
(zeros of $s \mapsto \zeta^{\circ n}(s) - s$) near the bulb
of $A_{\phi}$ containing $\rho_1$ for $2 \leq n \leq 5$.
As $n$ increases the pattern of
distribution of these zeros becomes
more and more obscure. 
The situation near $\rho_1$ appears to be typical of 
that near all nontrivial zeros of zeta iterates.
\newline  \newline
\begin{figure}[!htbp]
\centering
\includegraphics[scale=.4]{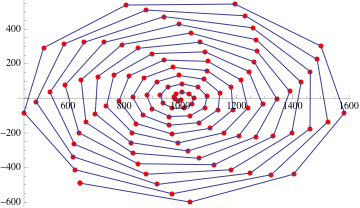}
\hskip .1in
\includegraphics[scale=.4]{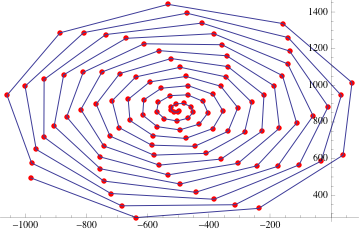}
\vskip .07in
\includegraphics[scale=.4]{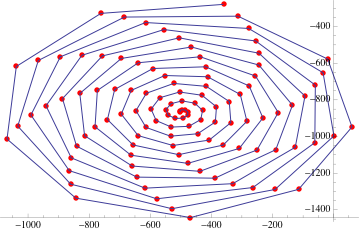}
\hskip .1in
\includegraphics[scale=.4]{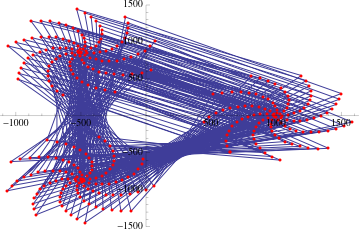}
\vskip .1in
\hskip 0in
\sc{fig. 5.9: the branch of a  $\rho_1$
 backward orbit induced by a zeta
$3$-cycle}
\end{figure}
 Figure 5.9 illustrates the branch 
 $B_{\rho_{_1},\Lambda}$ of
 the backward orbit of
  $\zeta^{\circ -}(\rho_1)$ induced by a $3$-cycle 
 $\Lambda = ( \lambda_1, \lambda_2, \lambda_3 )  $
 with $\lambda_1 \approx 3.95896 + 24.2362i$ . 
 The red vertices of a 
 given chord of the graph 
 represent points  in
 the branch. 
 Geometrically, Figure 5.9 is
doubly abstract: 
(1)  The spirals have been positioned so that their centers are
placed at 
$1000 + 0i, 1000(\cos(\frac{2\pi}{3})+ i \sin(\frac{2\pi}{3}))$
and $1000(\cos(\frac{4\pi}{3})+ i \sin(\frac{4\pi}{3}))$
 for the sake of legibility; and
(2) the spirals are everted: they have been 
re-scaled logarithmically, so that points closer to 
the center appear, in the figure, to be 
farther from the center.
The row 2, column 2 panel  shows vertices representing 
elements of $B_{\rho_{_1},\Lambda}$ 
and edges between vertex pairs
$(v, \zeta(v))$, while the other
panels depict the spirals separately;
these are portraits of 
$b_j (j = 1, 2, 3)$.
In these three figures, each vertex pair
$(v, \zeta^{\circ 3}(v))$  is connected by an edge.
\subsection{Angular distribution of branches along the spirals.}
A structural invariant, which seems to
determine the number of arms visible in
our plots of  branches of
$f^{\circ -}(z)$, is the  function 
 $\delta_{f,z, \psi}: k \mapsto \arg(a_k - \psi) - \arg(a_{k+1} - \psi)$,
where  the $a_k$'s are members of a particular 
branch  of $f^{\circ -}(z)$
 converging to $\psi$.
\newline \newline
The values of  $\delta_k$ in the case
 discussed in
section 5.2 correspond in the
present  notation to those of 
 $\delta_{\zeta,\,\rho_1,\, \psi_1}(k)$.
They appear to change very slowly
with $k$, and this behavior seems to be what gives rise to
the appearance of discrete arms in
plots of branches of $\zeta^{\circ -}(\rho_n)$.
\newline \newline
 Among 
the zeta fixed points very near 
trivial zeros, the only values
of 
\newline
$\delta_{\zeta,\rho_n, \psi_{_{-2n - 18}}}(k)$ that we see 
(Figure 5.7) are 
$\approx \pi$ and $\approx 2\pi$, distributed,
as we have noted above, according to the mod $4$
 residue classes of the zeros.
 \newline \newline
Because $\delta_{\zeta,\rho_n,\psi_{\rho_n}}(k)$ 
apparently converges
rapidly as $k \rightarrow \infty$ 
(as in Figure 5.5 where $n = 1$),
we take the  value of 
$\delta_{\zeta, \rho_n, \psi_{\rho_n}}(100)$ as a proxy for 
\newline
$\lim_{k\rightarrow \infty} \delta_{\zeta,\rho_n,\psi_{\rho_n}}(k)$.
\begin{figure}[!htbp]
\centering
\includegraphics[scale=.85]{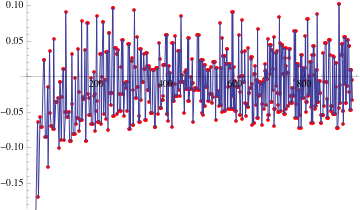}
\vskip .07in
\sc{fig. 5.10: $$\frac{\delta_{\zeta,\rho_n,\psi_{\rho_n}}(100)-\pi/2}{\pi/2},
1 \leq n \leq 600$$}
\end{figure}
Then our calculations
are consistent with the
proposition that
$\lim_{k\rightarrow \infty} \delta_{\zeta,\rho_n,\psi_{\rho_n}} \approx \pi/2$
(Figure 5.10.) Very small differences in this limit as $n$ varies appear
to determine very different shapes for the discrete arms visible
in our plots.
\newline \newline
We have observed in all of our
experiments that the visible  structure of
a branch of $f^{\circ -}(z)$ depends upon the fixed
point at its center  and not on $z$,  so
$\delta_{f,z, \psi}$ should  depend only upon $f$ and 
$\psi$. Contrary to the impression suggested
 by our notation, it should be independent of $z$,
but we cannot exclude the possibility
that there are counterexamples to this idea.
\newline \newline
\subsection{Logarithmic models of spirals interpolating branches of 
the backward orbit of zeta.}
The branches 
$B_{\rho, \psi} = (a_0 = \rho, a_1, a_2, ...)$ of
 $\zeta^{\circ -}(\rho)$ 
for nontrivial Riemann zeros $\rho$
converging to  zeta fixed points $\psi$
are interpolated by curves that
resemble logarithmic spirals.
We carried out experiments in which we 
looked for approximations
of these interpolating curves
by such  spirals. We chose 
 branches of the argument
function, varying with $k$ and evaluated at 
$a_k - \psi$, such that
an angle $\theta_k$ was assigned to 
$a_k - \psi$ which was the least such angle
 $ > c + \max_{j <k} \theta_j$
for  $c = 0$ or $1$.
 The angle
$\theta_0$ was the value 
of $\arg (a_0 - \psi)$
from the branch of argument
chosen automatically by \it Mathematica. \rm
 The choice of
$c$ was dictated by requiring
that $\theta_k$  act like
a winding number about $\psi$ 
evaluated at the points $a_k$. For $\psi$
near a trivial zero, $c = 1$ was chosen;
for $\psi$ near a nontrivial zero, $c = 0$.
\newline \newline
For $r_k =  |a_k -\psi|$
plots of the  sets of pairs 
$(\theta_k,\log r_k), k \geq 0$,
appear to lie on curves resembling
straight lines. For spirals centered at 
zeta fixed points $\psi_{\rho_n}$ near the $\rho_n$,
we approximated 
these lines 
using \it Mathematica\rm's FindFit
command.
\begin{figure}[!htbp]
\centering
\includegraphics[scale=.4]{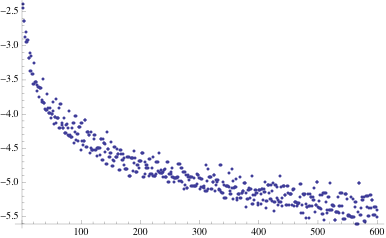}
\hskip .1in
\includegraphics[scale=.4]{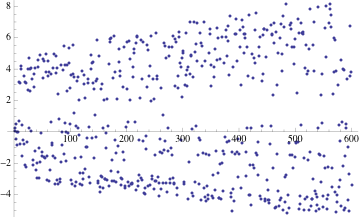}
\vskip .07in
\hskip 0in
\sc{fig. 5.11: slopes $m_n$ (left panel) and intercepts $b_n$
(right panel) in  log-linear models $\widetilde{r} = \exp (m_n\theta + b_n)$ of spirals interpolating
branches
of $\rho_n$ centered at $\psi_{\rho_n}$ with $r = |z  - \psi_{\rho_n}|$
and $\theta = \arg (z - \psi_{\rho_n})$}
\end{figure}
Figure 5.11 is a pair of plots of $m_n$ and $b_n$ against
$n$ for models $|z- \psi_{\rho_n}| = \exp (m_n \theta + b_n)$
 fitted to  branches of $\zeta^{\circ -}(\rho_n),
 1  \leq n \leq 600$. In particular,
 we write $\widetilde{r_k} = \exp (m_n \theta_k + b_n)$
 for our estimate of $r_k$. These seem to be
 first-order approximations to 
 genuine interpolating curves; the error-term
 will be
 discussed in the next section.
 \newline \newline
Investigating spirals centered at 
the zeta fixed points $\psi_{-2n}$ that
lie near the trivial zeros was carried out in
a different way. We were able to
 collect a substantial amount of
 data (meaning for the first $200$ members of the
 branches) for the
  $\psi_{\rho_n},  n \leq 600$ using
 $500$-digit precision. 
 On the other hand, even with 
 $1000$-digit precision, 
 we were able to collect data only
 on the first twenty elements
 of branches centered at the $\psi_{-2n}$ 
 for $n \leq 30$  before
 the use of the FindFit command
to get a linear model for the pairs
$(\theta_k, \log r_k)$
produced error messages
 from \it Mathematica. \rm
 \newline \newline
Fortunately, the branches spiraling about
the $\psi_{-2n}$ appear to
be better behaved than  the 
ones spiraling  about the $\psi_{\rho_n}$.
For each pair $(n, n^*)$
the branch $B_{\rho_n, \psi_{-2n^*}}$
is apparently 
interpolated both by a curve very nearly a 
logarithmic spiral and
by another curve which is very nearly
a straight line passing through 
the points $\rho_n$ and 
$\psi_{-2n^*}$.
As we will see in the next
section, the fit of the branches to
the straight line passing through
these two points
is so good
that we  use it as our second model,
together with the assumption
that the
 $\theta_k =  a\pi k \, +$ (a constant
depending only on $n$ and $n^*$),
with $a = 1$ or $2$ depending 
as we have explained
only on the parity of $n^*$. 
It was feasible to find 
linear models for
the maps $k \mapsto \log r_k.$
Combining the assumptions
about the $\theta_k$
with the linear models we 
construct for the $\log r_k$
gives logarithmic models for
the interpolating spirals.
We test these models in the next section.
\newline \newline
 We chose
 to examine the behavior of
 branches 
 $B_{\rho_n,\psi_{-2n-18}}$ 
 of $\zeta^{\circ -}(\rho_n)$
 because, among zeta fixed points close to the
 trivial zeros $-2n$, the greatest one
 (the one that lies rightmost
 along the real axis) is very close to
 $-20$.   
 \begin{figure}[!htbp]
 \begin{centering}
\includegraphics[scale=.4]{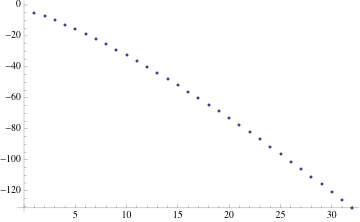}
\hskip .1in
\includegraphics[scale=.4]{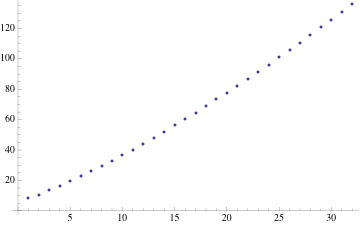}
\vskip .07in
\hskip 0in
\end{centering}
\sc{fig. 5.12: slopes $m_n$ (left panel) and intercepts $b_n$
(right panel) in  models $\widehat{r_k} = \exp (m_n k + b_n)$ of spirals interpolating
branches of $\rho_n$ centered at $\psi_{-2n-18}$}
\end{figure}
Figure 5.12 is a plot corresponding
 to Figure 5.11
 for the zeta fixed points $\psi_{-2n-18}$. 
 Here we have
 $B_{\rho_n,\psi_{-2n-18}} = (a_0 = \rho_n, a_1, a_2, ...)$.
 Writing $r_k = |\psi_{-2n-18} - a_k|$, we 
 took $m_n$ and $b_n$ to be, respectively, 
 the means of the slopes and intercepts of the chords connecting
 consecutive pairs $p_k = (k, \log |a_k - \psi_{-2n-18}|)$.
 Thus our model for $r_k = |a_k - \psi_{-2n-18}|$ is
 $\widehat{r_k} = \exp (m_n k + b_n)$.
 We discuss it further
 this in 
 the following section.
\section{\sc error terms}
We will take the phrase
``error term'' to encompass complex-valued
deviations from a given estimate as well as
their absolute values. 
Like the original
estimates, the curves followed by 
complex-valued deviations appear to have
the form of logarithmic spirals.
This raises the prospect of an infinite regress,
which might perhaps lead to an exact 
expression for the best interpolating curves,
but we have postponed any investigation of this
idea.
\subsection{Deviation of backward 
orbit branches from logarithmic spirals.}
\subsubsection{Branches converging to fixed points near non-trivial zeros.}
(This subsection provides some of our evidence
for Conjecture 1.) For nontrivial 
Riemann zeros $\rho_n$ and
the corresponding zeta fixed points $\psi_n$, we plotted the relative
complex-valued deviations  $d_{rel(\rho_n,\psi_n,a_k)}$.

\begin{figure}[!htbp]
\begin{centering}
\includegraphics[scale=.3]{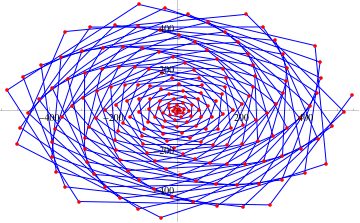} 
\hskip .0in
\includegraphics[scale=.3]{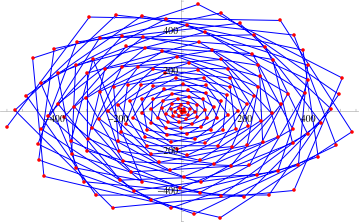} 
\hskip .0in
\includegraphics[scale=.3]{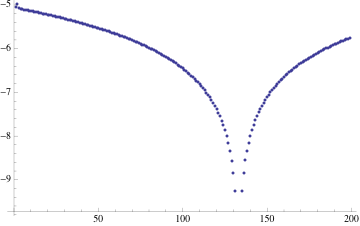} 
\vskip .7in
\includegraphics[scale=.3]{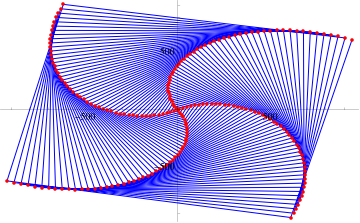}
\hskip .0in
\includegraphics[scale=.3]{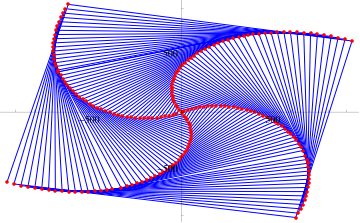}
\hskip .0in
\includegraphics[scale=.3]{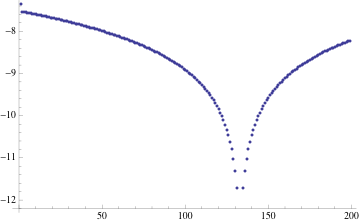}
\vskip .7in
\includegraphics[scale=.3]{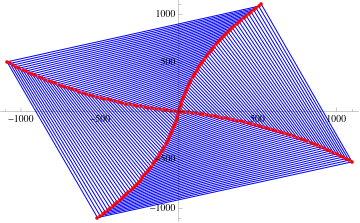}
\hskip .0in
\includegraphics[scale=.3]{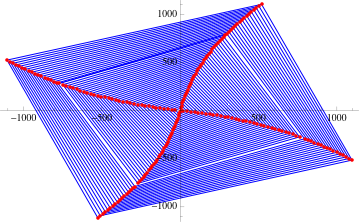}
\hskip .0in
\includegraphics[scale=.3]{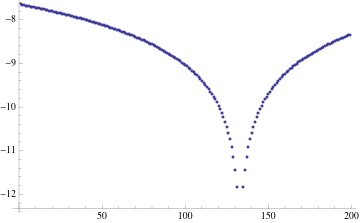}
\vskip .1in
\hskip .7in
\newline 
\vskip .1in
\hskip 0in
\end{centering}
\sc{fig. 6.1: $B_{\rho_n, \psi_n}, n = 1, 28, 48$: column 1: original branches;}
\newline
\sc{column 2: deviations $d_{rel(\rho_n,\psi_n,a_k)}$ of $B_{\rho_n, \psi_n}$ from logarithmic spirals;}
\newline
\sc{column 3: $\log d_{rel(\rho_n,\psi_n,a_k)}$ vs. $k$}
\end{figure}
 The runs depicted in Figure 6.1                
 portray both kinds of plots for 
 $\rho = \rho_n, \psi = \psi_{\rho_n}$ 
 for $n = 1, 28$, and $48$.
It seems noteworthy that
 the two kinds of plots
resemble each other so closely, but
 inspection demonstrates that they
are not identical.
\newline \newline
The values of $\log d_{rel(\rho_n,\psi_n,a_k)}$ for $n = 1, 28$, 
and $48$ are also plotted in
Figure 6.1 (column 3.)
The magnitude of the $d_{rel(\rho_n,\psi_n,a_k)}$
appears to decay exponentially
for  $k <$ roughly $130$; for larger $k$, 
however, the magnitude of the deviations
appears to grow exponentially without exceeding
$e^{-6}$ for $k \leq 200$. We think
that the shape of this curve,
which is typical, is an artificial
effect of the FindFit command on
a file of $200$ pieces of data: the fit is
best near the center of the data file.
\newline \newline
For a fixed point $\psi = \psi_{\rho_n}$ near a nontrivial 
zero $\rho = \rho_n$, 
we used the initial $200$ elements of
each branch $B_{\rho, \psi}$ as a proxy
for $B_{\rho, \psi}$
to study the 
relative deviations
of the curve interpolating it
from a logarithmic spiral.
Taking $\beta = 200$ for the moment, let us set
\newline \newline
$
max_{n, \beta} := \max_{1 \leq k \leq \beta} d_{rel(\rho_n,\psi_n,a_k)},
$
\newline \newline
$max_{n, \beta}^* :=  \sqrt{n/\log n} \times max_{n, \beta},$
\newline \newline
let $mean_{n, \beta}$ denote the mean of the  
$d_{rel(\rho_n,\psi_n,a_k)}, k = 1 ,2, ..., \beta$
and let 
\newline \newline
$mean_{n, \beta}^* :=  \sqrt{n/\log n} \times mean_{n, \beta}.$
\newline \newline
\begin{figure}[!htbp]
\centering
\includegraphics[scale=.4]{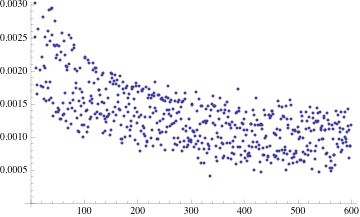}
\hskip .2in
\includegraphics[scale=.4]{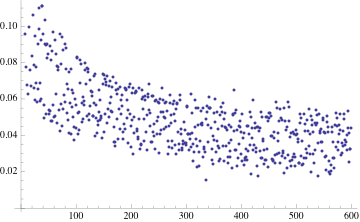}
\vskip .1in
\includegraphics[scale=.4]{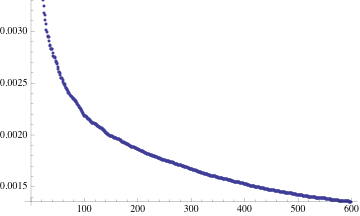}
\hskip .2in
\includegraphics[scale=.4]{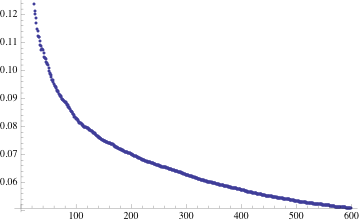}
\vskip .1in
\sc{fig. 6.2: row 1: $mean_{n, \beta}^*$ (left)
and  $max_{n, \beta}^*$ (right), $\beta = 200, 1 \leq n \leq 600$;}
\sc{\newline row 2: smoothings of the corresponding plots in row 1}
\end{figure}
(We exclude $k=0$ in these definitions;
thus the values of these four numbers
tell us nothing about the fit of $a_0 = \rho$
to the spiral in question.) 
\newline \newline 
The panel in row 1 column 1 of Figure 6.2 
is a plot of $600$ values of
$mean_{n, \beta}^*$.
 The panel in row 1 column 2  is a plot of 
$max_{n, \beta}^*$.
The plots in row 2 are  smoothings of the plots in 
row 1: for each $n$, they depict means of 
$mean_{j, \beta}^*$ and $max_{j, \beta}^*$ over the
range $1 \leq j \leq n$. 
\newline \newline
These plots are consistent with the proposition
that, for 
$\beta = 200$, $mean_{n, \beta}$ and $max_{n, \beta}$
both 
$= O \left ( \sqrt{\frac {\log n}{n}}\right )$
with both implied constants $<1$. More 
optimistically, perhaps, the plots are consistent
with the hypothesis that
 $mean_{n, \beta}$ and $max_{n, \beta}$
both 
$= o \left ( \sqrt{\frac {\log n}{n}}\right )$.
\newline \newline
We tested the same idea after replacing $\sqrt{\frac {\log n}{n}}$
with powers $\left(\frac {\log n}{n}\right)^{\epsilon}$
for $\frac 12 < \epsilon <1$. It seems possible that the
supremum of $\epsilon$ for which these statements 
might be true lies in the half-open interval $[.8, .9)$.
It also seems possible that this supremum is a 
decreasing function of $n$. We omit the 
relevant plots.
\subsubsection{Branches converging to fixed points 
near the trivial zeros.}
(This subsection provides 
some of our evidence for Conjecture 2.)
Let the branch $B_{\rho_n, \psi_{-2n-18}}$ of 
$\zeta^{\circ -}(\rho_n) = (a_0  = \rho_n, a_1, a_2, ....)$.
Before we test a logarithmic model for the decay of
$r_k = |a_k - \psi_{-2n - 18}|$,
we want to assess how well 
$B_{\rho_n, \psi_{-2n-18}}$ 
fits the straight line passing through 
$\rho_n$ and $\psi_{-2n-18}$.
We measured the vertical 
deviation of the $a_k \in B_{\rho_n, \psi_{-2n-18}}$
from the straight line passing through
both $\rho_n$ and $\psi_{-2n-18}$
as a fraction of the heights of the $a_k$.
We make the following definitions:
\newline \newline
$M_n$ and $B_n$ are the slope and intercept respectively
of the straight line passing through $\rho_n$ and
$\psi_{-2n  - 18}$
\newline \newline
and
\newline \newline
$d^{trivial}(n, k): =$
$$ \left|\frac{ \Im(a_k) - (M_n \Re(a_k) + B_n)}{\Im(a_k)}\right|,$$
\newline \newline 
$mean^{trivial}_{n, \beta} = $ the mean of the 
$d^{trivial}(n, k), 1 \leq k \leq \beta,$
\newline \newline
and
\newline \newline
$max^{trivial}_{n, \beta} = \max_{1 \leq k \leq \beta}d^{trivial}(n, k)$.
\newline \newline
\begin{figure}[!htbp]
\centering
\includegraphics[scale=.4]{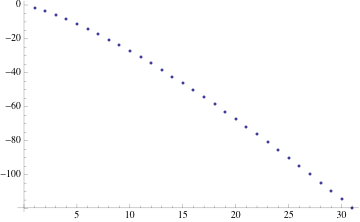}
\hskip .2in
\includegraphics[scale=.4]{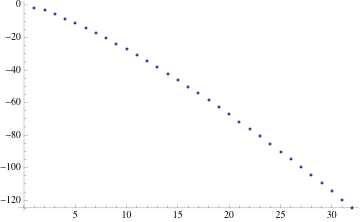}
\vskip .1in
\sc{fig. 6.3: $\log mean^{trivial}_{n, \beta}$ (left)
and  $\log max^{trivial}_{n, \beta}$ (right); $\beta = 20, 1 \leq n \leq 32$}
\end{figure}
The left panel of Figure 6.3 is a plot of 
$\log mean^{trivial}_{n, \beta}, 1 \leq n \leq 32$ and  $\beta = 20$.
The right panel is a corresponding  of $\log max^{trivial}_{n, \beta}$.
 Evidently, the $a_k$
lie near the specified lines, and agreement with the lines
improves rapidly as $n$ increases.
\newline \newline
Next we define  for $a_k \in B_{\rho_n, \psi_{-2n-18}}$
\newline \newline
$r_{n,k} = |a_k - \psi_{-2n-18}|$,
\newline \newline
for $k > 1, m_{n, k}$ = the slope of the chord connecting 
the ordered pairs $(k, \log r_{n, k})$ and $(k, \log r_{n, k-1})$ in
$\bf{R}\rm^2$,
\newline \newline
$m_{n, \beta} = $ mean  of the $m_{n, k}, 1 \leq k \leq \beta$,
\newline \newline
and we define a $y$-intercept function $b_{n, \beta}$ analogously.
The error functions are defined as follows:
$$d^{model}(n, k, \beta):=\left|\frac {\log r_{n,k} - (m_{n, \beta} k +
b_{n, \beta}) }{\log r_{n,k}}\right|,$$
\newline \newline
 $mean^{model}_{n, \beta} = $ the mean of the 
$d^{model}(n, k,  \beta), 1 \leq k \leq \beta,$
\newline \newline
and
\newline \newline
$max^{model}_{n, \beta} = \max_{1 \leq k \leq \beta} d^{model}(n, k,  \beta)$.
\newline \newline
\begin{figure}[!htbp]
\centering
\includegraphics[scale=.4]{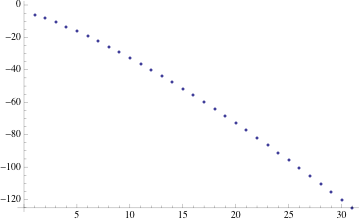}
\hskip .2in
\includegraphics[scale=.4]{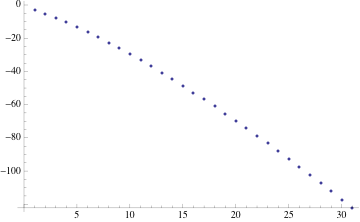}
\vskip .1in
\sc{fig. 6.4: $\log mean^{model}_{n, \beta}$ (left)
and  
$\log max^{model}_{n, \beta}$ (right); $\beta = 20, 1 \leq n \leq 32$}
\end{figure}
The left panel of Figure 6.4 is a plot of 
$\log mean^{model}_{n, \beta}, 1 \leq n \leq 32$ and  $\beta = 20$.
The right panel is a corresponding  of $\log max^{model}_{n, \beta}$.
Once more the fit is good and improves rapidly as $n$ increases.
\subsubsection{Deviation of the Riemann zeros from
fitted logarithmic spirals.}
The plots in Figure 6.5  display
numerical information that support 
Conjecture 4 and the scenario described in section 6.1.3.
They indicate the 
possibility that as $n \to \infty$, the nontrivial Riemann zeros
$\rho_n$
become more well fitted to the logarithmic spirals
we have in turn 
fitted to the  branches $B_{\rho_n,\psi_n}$.
The left panel  is a plot of $\log d_{rel}(\rho_n,\psi_n,0)$ against 
$n, 1 \leq n \leq 600$;
The right panel plots $D_{rel}(N)$, the log of the running mean
of the  $d_{rel}(\rho_n,\psi_n,0)$, as defined in
Conjecture 4, 
$1 \leq N \leq 600$.
Figure 6.6 shows the corresponding plots 
for  $d_{rel}(\rho_n,\psi_{n+1},0)$. 
\begin{figure}[!htbp]
\centering
\includegraphics[scale=.4]{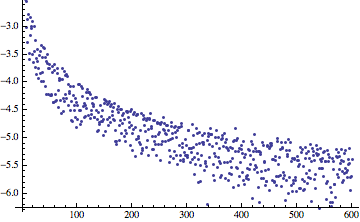}
\hskip .2in
\includegraphics[scale=.4]{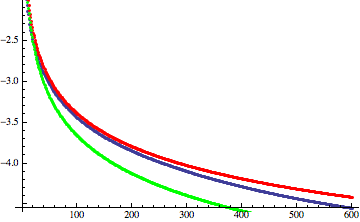}
\vskip .1in
\sc{fig. 6.5: left, $\log d_{rel}(\rho_n,\psi_n,0), 1 \leq n \leq 600$};
\newline
\sc{right, $\log \frac 1N \sum_{n=1}^N d_{rel}(\rho_n,\psi_n, 0), 1 \leq N \leq 600$};
\newline
\sc{green curve: $-(\log N)^{.8}$, red curve :$-(\log N)^{.85}$}
\end{figure}
\begin{figure}[!htbp]
\centering
\includegraphics[scale=.4]{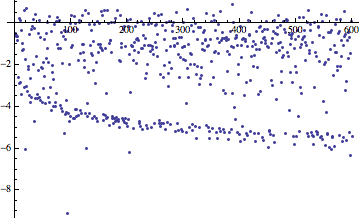}
\hskip .2in
\includegraphics[scale=.4]{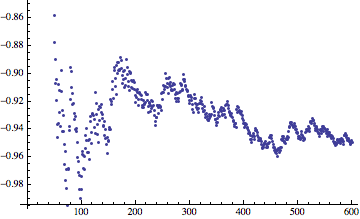}
\vskip .1in
\sc{fig. 6.6: left, $\log d_{rel}(\rho_n,\psi_{n+1},0), 1 \leq n \leq 600$};
\newline
\sc{right, $\log \frac 1N \sum_{n=1}^N d_{rel}(\rho_n,\psi_{n+1}, 0), 1 \leq N \leq 600$};
\end{figure}
The next two plots treat absolute deviations $d_{abs}$, which 
suggest the  narrowing of the widths of the ``error bands''
mentioned in the  speculations we ventured in the Introduction.
Figures 6.7 corresponds to the plots of relative deviations in Figure 6.5.
\begin{figure}[!htbp]
\centering
\includegraphics[scale=.4]{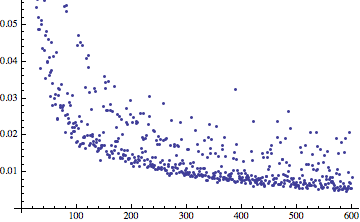}
\hskip .2in
\includegraphics[scale=.4]{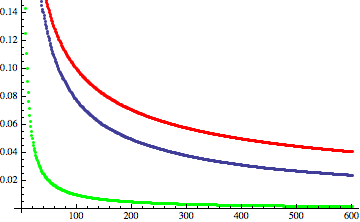}
\vskip .1in
\sc{fig. 6.7: left, $\log d_{abs}(\rho_n,\psi_n,0), 1 \leq n \leq 600$};
\newline
\sc{right, $\log \frac 1N \sum_{n=1}^N d_{abs}(\rho_n,\psi_n, 0), 1 \leq N \leq 600$};
\newline
\sc{green curve: $\frac 1N$, red curve :$\sqrt{\frac1N}$}
\end{figure}
Figure  6.8 displays information on absolute	
deviations for fixed points $\psi_{-2n}$ nearest to the
trivial zeros $-2n$, with spirals terminating at $\rho_n$ (left panel)
and $\rho_{n+1}$ (right panel).
\begin{figure}[!htbp]
\centering
\includegraphics[scale=.4]{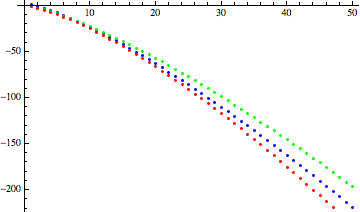}
\hskip .2in
\includegraphics[scale=.4]{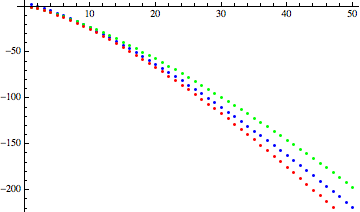}
\vskip .1in
\sc{fig. 6.8: left, $\log d_{abs}(\rho_n,\psi_{-2n-18},0), 1 \leq n \leq 600$;}
\newline
\sc{right, $\log d_{abs}(\rho_{n+1},\psi_{-2n-18},0), 1 \leq n \leq 600$;}
\newline
\sc{both green curves: $-n^{1.35}$, both red curves : $-n^{1.4}$}
\end{figure}
\subsection{Deviation from rotational invariance.}
It did not seem 
plausible to us that
something special about Riemann zeros $\rho$
should force the 
branches $\zeta^{\circ -}(z), z = \rho$
 in particular, to be attracted
to  repelling fixed points  $\psi$
along logarithmic spirals.
As we remarked in the introduction, dynamical
systems theory leads one to expect
that the $\psi$ should attract
 all  of the nearby branches $\zeta^{\circ -}(z)$
whether or not $\zeta(z) = 0$,
and (one speculated) probably along 
roughly similar curves.
Logarithmic spirals appear 
in fluid mechanics 
(see, \it e.g.\rm, \cite{Ma}, v.2, pp. 186-188
or \cite{SpAk}, p. 358.)
By analogy with the streamlines of
a vortex  in a fluid,
we  speculated that
 the
existence of spiral curves connecting
zeros 
 of zeta
to  repelling zeta fixed points might
be a consequence of 
a scenario in which 
there is
an infinite family of such 
spirals  related by rotations around
the fixed point. 
By this we mean
a family of
spirals $s$ parameterized by real numbers $x$, 
varying continuously with $x$
in a sense  made explicit by condition (1) below,
such that if $s_x$ and $s_{x + \theta}$
are two such spirals with common center a
zeta fixed point
$\psi$, then 
\newline \newline
(1)
$ s_{x + \theta} - \psi = e^{i \theta} ( s_x -\psi)$
\newline \newline
(condition (1) being an equation of homotheties) and
\newline \newline
 (2)
$z \in s_x \Rightarrow$  
(i) $\zeta(z) \in s_x$ and 
 (ii) there exists a
branch $B_z \subset s_x$ of $\zeta^{\circ -}(z)$
such  that $\lim B_z = \psi$.
\newline \newline
In  this scenario, the spiral curves would be
congruent in the sense of
Euclidean geometry and  exactly one
spiral
 would intersect the critical line
at each $\rho_n$
without appealing to special properties 
of the zeros. We have verified the existence of
spiral branches $B_z$ 
of $\zeta^{\circ -}(z)$ for various $z$  on the 
critical line other
than Riemann zeros  without meeting a counterexample.
Like the spirals we have already described,
they are approximately logarithmic;
we omit the relevant plots.
\newline \newline
Condition (1) imposes rotational invariance
on the $s_x$. This suggests the  
possibility that the branches of $\zeta^{\circ -}$
interpolated by them enjoy the
same property. Suppose 
$u$ and $\zeta(u)$ lie on $s_x$ with center
 $\psi$ and
let $R_{\theta,\psi}(z) := e^{i \theta}(z - \psi) + \psi$ 
be the  function that takes $z$ to its image under
rotation by an angle $\theta$ around $\psi$.
Under perfect rotational invariance, not
only of the spirals $s_x$ but of the branches
of $\zeta^{\circ -}(z)$ for particular $z$
that they interpolate, the numbers 
$R_{\theta,\psi}(\zeta(u)) - \zeta(R_{\theta,\psi}(u))$
must vanish. Therefore we studied the 
$R_{\theta,\psi}(\zeta(u)) - \zeta(R_{\theta,\psi}(u))$.
We restricted ourselves to $u \in \zeta^{\circ -}(z)$
for various $z$,
not only because these are the main objects of
interest, but because our only 
reliable information
about the spirals
comes from their interpolation
of the branches $B_{\psi} =
(a_0, a_1, a_2, ...)$
 of $\zeta^{\circ -}(z)$,
and so our only useful  candidates for points in $s$
are the members of such branches. 
\newline \newline
Logarithmically scaled plots (which we omit) of the
discrepancies (say)  \newline
$R_{\theta,\psi}(\zeta(a_n)) - \zeta(R_{\theta,\psi}(a_n))$
indicate that these numbers decay in modulus
exponentially and rotate
around the origin in a  nearly linear fashion with $n$.
In other words, they themselves describe curves
that are approximated by logarithmic
spirals. 
\section{\sc appendix: the figures}
\subsection{Figure 1.1}
Figure 1.1 a 120 
by 120 square with center $1 + 0i$.
\subsection{Figure 2.1}
Figure 2.1 depicts a  $6$ by $6$ square with center zero. 
\subsection{Figure 3.1}
 Figure 3.1 shows an 8 by 8 square
 with center $-5 + 9.5i$.
\subsection{Figure 3.2}
Figure 3.2 shows  a
120 by 120 square with center zero.
\subsection{Figure 3.3}
Figure 3.3 depicts 
a 60 by 60 square
with center zero.
\subsection{Figure 3.4}
The upper left panel of Figure 3.4 
shows a $2.4 \times 10^{-5}$
by  $2.4 \times 10^{-4}$ square
 centered at $-28$.
The upper right panel shows a $2.4 \times 10^{-4}$
by  $2.4 \times 10^{-4}$ square
 centered at $-26$. The lower left
 panel shows a $.004$ by $.004$ square 
 centered at $-24$. The lower right depicts a 
 $.07$ by $.07$ square centered at $-22$. 
\subsection{Figure 3.5}
 The panel in row 1 column 1
 of Figure 3.5 depicts  a $10$ by $10$
square centered at $\rho_1$.
The other panels show a square
with side length $.006$ 
and center $\rho_1 + 4.1215 - .4015i \approx 4.6215 + 13.7332 i$.
\subsection{Figure 3.6}
 Each panel of Figure 3.6 is a 30 by 30 square with center $ = -5$.
\subsection{Figure 4.1}
 The squares depicted in Figure 4.1
 have side length
 $.2, .02, .002,$ and $.0002$ 
 in rows 1, 2, 3 and 4, respectively.
 The center of each square is 
 $\psi_{_{\rho_1}} \approx -2.3859 + 16.271i$.
\subsection{Figure 5.1}
 All the panels of Figure 5.1
 depict $A_{\phi}$ in
  $2$ by $2$ squares. In rows 1 - 4,
 the centers are $\psi_{_{\rho_1}}$ -  $\psi_{_{\rho_4}}$, 
 respectively, where 
 $\psi_{_{\rho_2}} \approx -2.0369 + 21.9931i$,
  $\psi_{_{\rho_3}} \approx -1.6935 + 26.5283i$,
  and
   $\psi_{_{\rho_4}} \approx -1.7496 + 30.8158i$.
\subsection{Figure 5.2}
 In Figure 5.2, column 1 depicts branches of
 $\zeta^{\circ -}(\rho_n), 1 \leq n \leq 5$,
 centered at a zeta fixed point 
 $\approx -14.613 + 3.108 i$;
 column 2 depicts branches of 
  $\zeta^{\circ -}(\rho_n), 1 \leq n \leq 5$,
 centered at a zeta fixed point
 $\approx -5.279 + 8.803i$.
\subsection{Figure 5.3}
 In Figure 5.3 (referring to the caption), 
 the value of $n$ in row $a$, column
$b$ is $n = 2a + b - 2$. 
\subsection{Figure 5.8}
 Figure 5.8 depicts four views of
 a $12$ by $12$ square
 with center $\rho_1$.
\subsection{Figure 5.9}
 Figure 5.9 depicts the branch 
 $B_{\rho_{_1},\Lambda}$ of
 the backward orbit of
  $\zeta^{\circ -}(\rho_1)$ induced by a $3$-cycle 
 $\Lambda = ( \lambda_1, \lambda_2, \lambda_3 )  $
 with $\lambda_1 \approx 3.95896 + 24.2362i$. 
\newpage
\bibliography{bibtexcite}
$\mbox{barrybrent@member.ams.org}$
\end{document}